\documentclass[11pt]{amsart}

\usepackage{graphicx}
\usepackage{amscd}

\usepackage{fullpage}

\newtheorem{theorem}{Theorem}
\theoremstyle{plain}

\newtheorem{corollary}{Corollary}

\newtheorem{definition}{Definition}
\newtheorem{example}{Example}

\newtheorem{lemma}{Lemma}

\newtheorem{proposition}{Proposition}

\usepackage{amsmath,amssymb}  

\begin{document}

\title{
Appell Polynomials and Their Zero Attractors
}

\author{Robert  P. Boyer}
\address{  {\rm R. Boyer:}  Department of Mathematics, Drexel University,
Philadelphia, PA}
\email{rboyer@drexel.edu}
\author{William M. Y. Goh}
\address{ {\rm  W. M. Y. Goh}: Department of Statistics and Finance\\
University of Science and Technology of China (USTC)\\
Hefei 230026,China}
\email{ wmygoh@hotmail.com}
\date{ \today }

\begin{center}
{
\sc
\large

}
\end{center}

\begin{abstract}
A polynomial family $\{ p_n(x)\}$ is Appell if it is given  by 
$\frac{ e^{xt}}{g(t)} = \sum_{n=0}^\infty p_n(x) t^n$ or,
equivalently, $p_n'(x) = p_{n-1}(x)$. If $g(t)$ is an
entire function, $g(0)\neq 0$, with at least one zero,
 the asymptotics of linearly scaled polynomials $\{p_n(nx)\}$ are described by means of
 finitely zeros of $g$, including those of minimal modulus.
 As a consequence,
we determine the limiting behavior of their zeros as well
as their density. The techniques and results extend our earlier work
on Euler polynomials.
 \end{abstract}

\maketitle
\tableofcontents

\section{Introduction}

Let $g(t)$ be an entire function such that $g(0) \neq 0$.

\begin{definition}
The Appell polynomials $\{p_n(x)\}$ associated with generating function $g(t)$ are given by
\begin{equation}\label{eq:gen_fct}
\frac{ e^{xt}}{g(t)} = \sum_{n=0}^\infty p_n(x) t^n.
\end{equation}
\end{definition}

Some important examples are: the Taylor polynomials of $e^x$, with $g(t)=1-t$;  the Euler polynomials, with $g(t)=(e^t+1)/2$;
and Bernoulli polynomials, with $g(t)= (e^t-1)/t$;  and their higher order analogues.

The asymptotics and limiting behavior of the zeros of these families have been investigated
by many people; for example, \cite{g-b}, \cite{szego}, and so on.

In this paper, we obtained the asymptotics and the limiting behavior of the zeros
for all Appell families provided the generating function $g(t)$ satisfies one further condition:
that $g$   must have at least one zero.
  We use the ideas in our earlier paper \cite{g-b}; furthermore,
we simultaneously simplify and generalize some of the techniques there.

We found that the asymptotics in the general case are build from the basic example
$g(t)=1-t$ which coincides with the classical work of Szeg\"o on the Taylor polynomials
of the exponential function.  In our paper \cite{g-b}, we found that the asymptotics
for the Euler and Bernoulli polynomials are controlled by certain roots of $g(t)$, the
ones of minimal modulus.  In the general situation, as expected, the minimal modulus
roots of $g(t)$ are needed to describe the asymptotics but there may be finitely many  other roots needed
to determine the asymptotics.
These additional roots are determined through a geometric condition described in the terms of
rotated and scaled versions of the Szeg\"o curve: $ | xe^{1-x}|=1$, $|x|\leq 1$, $x \in {\mathbb C}$
(see Figure \ref{fig:szego_curve}).

We frequently use the following notations.
Let $Z(g)$ denote the set of all zeros of $g$ and let $r_0 < r_1 < r_2 < \dots$ denote
the distinct moduli of these zeros in increasing order. 

Recall that if $K_1$ and $K_2$ are two non-empty compact subsets of $\mathbb C$,
then their Hausdorff distance is the larger of $\sup\{ d(x,K_1) : x \in K_2\}$ and
$\sup\{ d(x,K_2) : x \in K_1\}$.

\begin{definition}
For a family $\{ q_n(x) \}$ of polynomials whose degrees are increasing to infinity,
their zero attractor is the limit of their set of zeros $Z(q_n)$ in the {\sl Hausdorff metric} on
the space of all non-empty compact subsets of the complex plane $\mathbb C$.
\end{definition}

In the appendix, we discuss how the zero attractor is found in terms of the
limsup and liminf of the zero sets.

\begin{figure}[h]
\includegraphics[height=3.0in, width=3.0in]{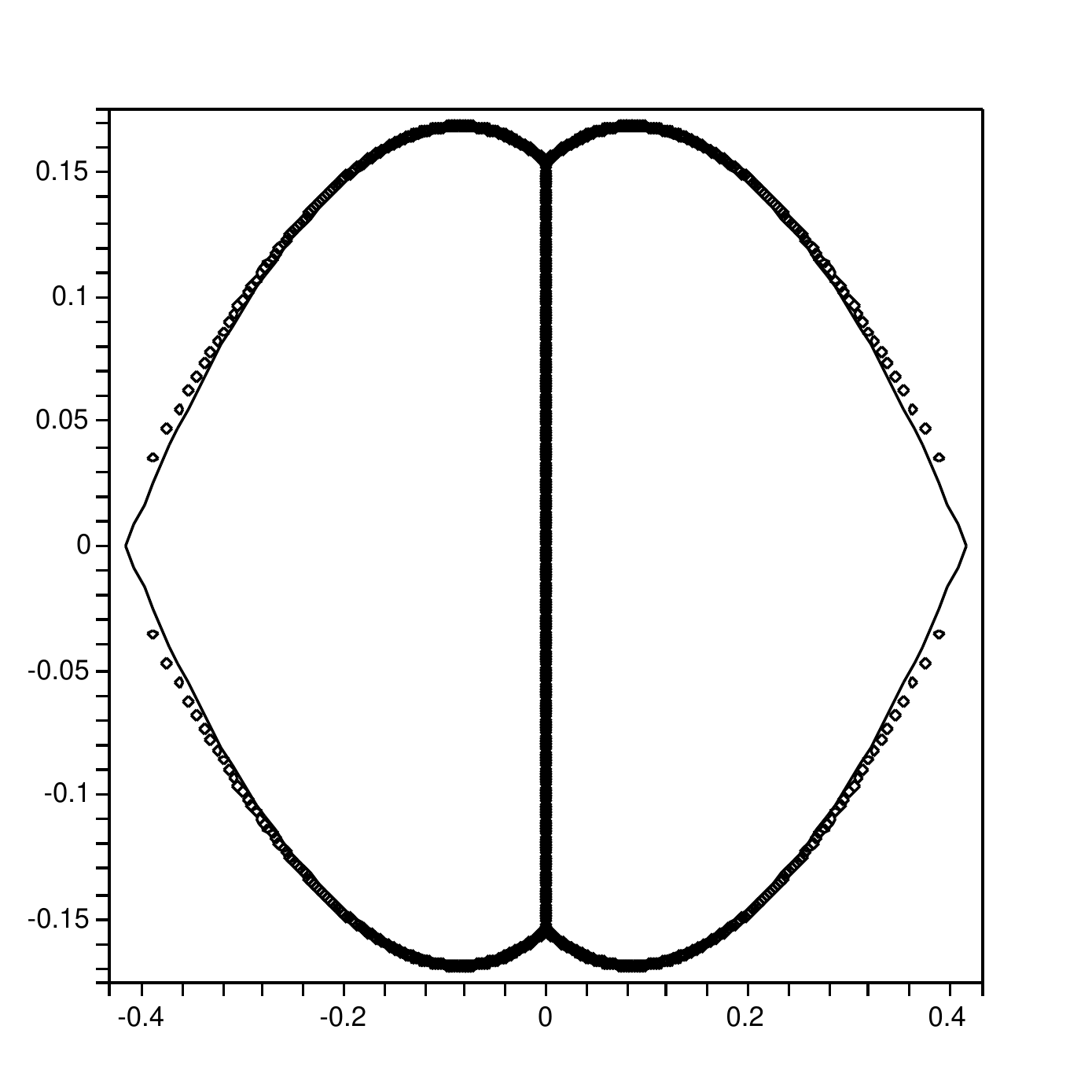} 
\caption{Zeros for degree 1000 polynomial, with generating function $g(t)=J_0(t)$ }
\label{fig:bessel_attractor}
\end{figure}

There is related work on the asymptotics and zeros of the Taylor polynomials
for linear combinations of exponentials $\sum c_j e^{\lambda_j x}$ where the 
parameters $\lambda_j$ satisfy a geometric constraint \cite{bleher}. The techniques of proof
are very different from our approach.

\section{The Generalized Szeg\"o Approximations}  % section 2

It is convenient to collect together several results from \cite{g-b} and some
extensions of them concerning the asymptotics of $S_n( x) = \sum_{k=0}^n x^k/k!$.
The domains of where their asymptotics hold are critical in understanding the
behavior for the Appell polynomials.

\begin{proposition}
{\sc (Left-Half Plane)}
\label{prop:4.1}
%\label{approximation-1} 
Let  $1/3<\alpha <1/2$ and $1 \leq j$.
On any compact subset $K$ of $\{w : \Re w <1 \}$, we have 
\begin{enumerate}
\item
$\displaystyle
\frac{S_{n-1}(nw)}{e^{nw}}=1-\frac{(we^{1-w})^{n}}
{\sqrt{2\pi n}(1-w)} \left(1+O(n^{1-3\alpha }) \right)$,
\item
$\displaystyle
D_{w}^{j-1}(w^{-n}S_{n-1}(nw))=D_{w}^{j-1}(w^{-n}e^{nz})
-
\frac{(j-1)!}{\sqrt{2\pi n}}\frac{e^{n}}{(1-w)^{j}} \left(1+O(n^{1-3\alpha }) \right)$,
\\
where 
the big O constant holds uniformly for $x \in K$.
\end{enumerate}
\end{proposition}

The proof of part (a) is in \cite{g-b}. Part (b) follows from an application
of the saddle point method.

The following Proposition is  also from \cite{g-b}:

\begin{proposition}
{\sc (Outside Disk)}
\label{prop:4.2}
%\label{approximation-2}
Let $S$ be a subset contained in $\left| w\right| >1$ with distance $\delta>0$ from
the unit circle,
and let $\alpha $ be chosen so $1/3<\alpha <1/2$. Then 
\begin{equation*}
\frac{S_{n-1}(nw)}{e^{nw}}=\frac{(we^{1-w})^{n}}{\sqrt{2\pi n}(w-1)}
\left(1+O(n^{1-3\alpha }) \right),
\end{equation*}
where the  big $O$ term holds uniformly for $w\in S$.
\end{proposition}

\begin{proposition}
{\sc (Evaluations of Integrals)}
\label{prop:4.3}
%\label{sn}
If $\epsilon <\left| w\right|$ and  $j\geq 1$, then we have
\begin{enumerate}
\item 
$\displaystyle
\frac{1}{2\pi i}\oint_{\left| t\right| =\epsilon }
\left(
\frac{e^{xt}}{t}\right)^{n}  \frac{1}{t-w}  \,dt= -w^{-n}S_{n-1}(wxn)$.
\\
\item
$\displaystyle
\frac{1}{2\pi i}\oint_{\left| t\right| =\epsilon }
\left(\frac{e^{xt}}{t}\right)^{n}\frac{1}{(t-w)^{j}} \,dt=\frac{-1}{(j-1)!}
D_{w}^{j-1}(w^{-n}S_{n-1}(wxn))$,\\
where $D_{w}$ is the differentiation operator $\frac{d}{dw}$.
\end{enumerate}
\end{proposition}
\begin{proof} %  proof 9  proposition
 (a) By expanding ${1}/{(t-z)}$ into an infinite geometric series
and performing a term-by-term integration, we obtain
\begin{eqnarray*}
\frac{1}{2\pi i}\oint_{\left| t\right| =
\varepsilon }\left(\frac{e^{xt}}{t}\right)^{n}
\frac{1}{t-z} \,dt
&=&
\frac{-1}{z2\pi i}\oint_{\left| t\right| =\varepsilon }
\left(
\frac{e^{xt}}{t}
\right)^{n}\frac{1}{1-\frac{t}{z}} \, dt
\\
&=&
\frac{-1}{z2\pi i}\oint_{\left| t\right| =\varepsilon }
\left(\frac{e^{xt}}{t}\right)^{n} \,
\left(\sum_{m\geq 0} \left(\frac{t}{z} \right)^{m} \right) \, dt.
\end{eqnarray*}
By the Cauchy integral theorem the terms correspond to $m\geq n$ vanish. 
Hence
\begin{eqnarray*}
\frac{1}{2\pi i}\oint_{\left| t\right| =
\varepsilon }\left(\frac{e^{xt}}{t}\right)^{n}
\frac{1}{t-z}\,dt
&=&
\frac{-1}{z}\sum_{n-1\geq m\geq 0}\frac{1}{z^{m}}
\left(\frac{1}{2\pi i} \,
\oint_{\left| t\right| =\varepsilon }e^{xtn}t^{-n+m} \, dt \right)
\\
&=&
\frac{-1}{z}\sum_{n-1\geq m\geq 0}\frac{1}{z^{m}}
\frac{(xn)^{n-m-1}}{
(n-m-1)!}
\\
&=&
\frac{-1}{z}z^{-n+1}\sum_{n-1\geq m\geq 0}\frac{(xnz)^{n-m-1}}{
(n-m-1)!} \\
&=&
-z^{-n}\sum_{n-1\geq j\geq 0}\frac{(xnz)^{j}}{j!}=(-1)z^{-n}
S_{n-1}(zxn).
\end{eqnarray*}
Part (b) follows from differentiating (a) $j-1$ times with respect to
$z$.
\end{proof} %proposition

\section{Asymptotics Outside the Disk $D(0;1/r_0)$}  % section 3

\begin{theorem}
\label{thm:3.1}
Let $K$ be any compact subset in the annulus $A(1/r_0; \infty)$.  
We have
\[
\frac{ p_n(nx)}{ (xe)^n/\sqrt{2\pi n}} = \frac{1}{g(1/x)} \left( 1 + O(1/n) \right).
\]
holds uniformly for $x\in K$.
\end{theorem}

\begin{proof} 
We shall find an asymptotic approximation to $p_{n}(nx)$ in the region 
$A(1/r_0; \infty)=\left\{x:\left| x\right| >\frac{1}{r_0} \right\}$. Use the generating relation
equation (\ref{eq:gen_fct}) to get
\begin{equation*}
p_{n}(x)=
\frac{1}{2\pi i}\oint_{\left| t\right| =
\epsilon }\frac{e^{xt}}{g(t)t^{n+1}} \, dt,
\end{equation*}
where $\epsilon <r_0$. Since both sides of the above equation are
entire functions of $x$, by analytic continuation this representation for 
$p_{n}(x)$ is valid for all $x\in C$. Hence we can replace $x$ by $nx$ to get
\begin{equation}
p_{n}(nx)=\frac{1}{2\pi i}\oint_{\left| t\right| =\epsilon }
\left(\frac{e^{xt}}{t}\right)^{n} \, \frac{dt}{tg(t)}.  \label{eq:pnnx-1}
\end{equation}
The above expression is valid for $0<\epsilon <r_{0}$ and is the starting
point of the analysis in the sequel. 

Let $K$ be an arbitrary compact
subset $\subseteq \{x:\left| x\right| >\frac{1}{r_{0}}\}$ and let $x\in K$.
We can certainly choose $\epsilon $ small enough so that for all $x\in K$,
$\left| \epsilon x\right| <1$. By a change of variables, we get
\begin{equation*}
p_{n}(nx)=\frac{x^{n}}{2\pi i}\oint_{\left| t\right| =\epsilon \left|
x\right| }
\left(\frac{e^{t}}{t}\right)^{n}\frac{dt}{tg(t/x)}.
\end{equation*}
Observe that the zeros of $g(t/x)$ have the form  $ax$  where $a \in Z(g)$.
Moreover, they must lie outside the closed unit disk since $|x|>1/r_0$,
so we can deform the integration path from the circle with
radius $\epsilon \left| x\right| $ to the unit circumference. Thus
\begin{eqnarray*}
p_{n}(nx)
&=&
\frac{x^{n}}{2\pi i}\oint_{\left| t\right| =1}
\left(\frac{e^{t}}{t}\right)^{n}
\frac{dt}{tg(t/x)}\\
&=&
\frac{x^{n}}{2\pi i}\oint_{\left| t\right| =1}
e^{n(t-\ln t)}\frac{dt}{tg(t/x)}.
\end{eqnarray*}
It can be easily seen that $t=1$ is the saddle point of the integral and the
classical saddle point method is applicable here  \cite{copson}. Hence
\begin{equation*}
p_{n}(nx)=\frac{(ex)^{n}}{\sqrt{2\pi n}g(\frac{1}{x})}\left(1+O(\frac{1}{n})\right),
\end{equation*}
where the implied $O$ constant holds uniformly for $x\in K$. 
\end{proof}  

The last equation can be written as
\begin{equation*}
\frac{p_{n}(nx)}{(ex)^{n}/\sqrt{2\pi n}}=\frac{1}{g(\frac{1}{x})}
\left(1+O(\frac{1}{n}) \right), \quad |x| > 1/r_0.
\end{equation*}

We have the:

\begin{corollary}\label{cor:3.2}
(a) On the complement of the disk $D(0;1/r_0)$,
$\displaystyle \lim_{n\to \infty} \frac{1}{n}
\ln \left| \frac{p_{n}(nx)}{(ex)^{n}/\sqrt{2\pi n}} \right| =0$
where the limit holds uniformly on compact subsets.
\\
(b)
The zero attractor must be contained in the closed disk
$\overline{D}(0; 1/r_0)$.
\end{corollary}

Note that part (b) follows easily from (a) since $g(x)$ never vanishes on the disk $D(0;1/r_0)$.

\section{Asymptotics on the Basic Regions $R_\ell$}   % section  4

Let $r_0, r_1, \dots$ denote the distinct moduli of the zeros $a$ of the generating function
$g$. Fix an integer $\ell$.
We fix  $\rho>0$ so it is not equal to any zero modulus $\{r_0,r_1,\dots\}$.
For each zero $a\in Z(g)$ with $|a|=r_\ell$,
 we consider the circle $|x|=1/|a|$ and the  disk $D(1/a, \delta_a)$.
 
Now the  tangent line $T_a$ to the circle $|x|=1/|a|$ at the point $1/a$
determines  the half-plane  $H_a$, which contains $0$; that is, $\Re(ax)<1$.
We choose $\epsilon_\ell>0$
 to be less than the distance from the portion of the tangent line $T_a$ that lies
outside the disk $D(1/a;\delta_a)$ to the circle $|x|=1/|a|$ for any $|a| = r_{\ell+1}$;
that is, $\epsilon_\ell < \sqrt{ 1/ r_\ell^2 + \delta_a^2} - 1/r_\ell$.
Finally, we make the requirement  the  disks $D(1/a;\delta_a)$ be mutually disjoint
for all $a \in Z(g)$ with $|a| < \rho$.

\begin{definition}
With these conventions, the region $R_\ell$ is described in terms of the half-planes $H_a$
and disks as
\begin{equation}\label{def:R_ell}
R_\ell =
\bigcap 
\left\{
H_a \setminus D(\tfrac{1}{a}; \delta_a) : |a| = r_\ell
\right\} 
\setminus    D(   0; \tfrac{1}{r_{\ell+1}} + \epsilon_{\ell+1} )
\end{equation}
\end{definition}

We note that the regions $R_\ell$ are not disjoint; in fact, by construction,
its inner boundary which consists  of 
the portion of the circle $|x| = \tfrac{1}{r_{\ell+1}}+ \epsilon$ that lie outside
the disks $D(1/a; \delta_a)$, $|a|=r_{\ell+1}$, actually lies inside the region $R_{\ell+1}$.

It is convenient to introduce a region that contains all of the $R_\ell$'s.

\begin{definition}
Let $R_\rho$ be the domain given as
\begin{equation}
R_\rho
=
\bigcap \left\{
H_a : a \in Z(g), |a|=r_0
\right\}  \setminus
\left[
\bigcup \left\{ D(1/a; \delta_a) : a\in Z(g), |a| < \rho \right\} \cup D(0; 1 / \rho)
\right].
\end{equation}
\end{definition}

Note the order of dependence:  first we can given the cut-off modulus $\rho>0$ for the moduli
of the zeros; next, $\delta_a>0$ for each $a \in Z(g)$ is given and is a function
of $\rho$ [see later section], then finally, $\epsilon_{\ell}$ is 
determined relative to each zero moduli $r_\ell$ which is a function of $\delta_a$.

 For any $a \in Z(g)$
with $r_0 \leq |a| < \rho$, let $s_a(t)$ be the singular part of
\begin{equation*}
\frac{1}{ t g(t)}
\end{equation*}
at its pole $a$. Next we set $g_1(t)$ to be
\begin{equation}\label{eq:def_g_1}
g_1(t) = \frac{1}{ t g(t)} 
-
\sum \{ s_a(t) : a \in Z(g), r_0 \leq |a| < \rho\}
\end{equation}
we see that $g_1(t)$ is analytic for $| t| < \rho$.

We develop the asymptotics for $\{ p_n(nx)\}$ on the regions $R_\ell$ where
$r_0 \leq r_\ell  <  \rho$.
Now
we saw already that we can write $p_n(nx)$ as
\[
p_n(nx) =
\frac{1}{2 \pi i} \int_{ |t| = \epsilon} \left( \frac{ e^{xt}}{t} \right)^n g_1(t) \, dt
+
\frac{1}{2 \pi i} \int_{ |t| = \epsilon} \left( \frac{ e^{xt}}{t} \right)^n s(t) \, dt,
\]
where $s(t) = \sum \{ s_a(t) : a \in Z(g), r_0 \leq |a| \leq \rho\}$.

\begin{lemma} \label{lemma:4.1}
With $g_1(t)$ given above in equation \ref{eq:def_g_1}, we have
\[
\frac{1}{2 \pi i}
\int_{ |t| = \epsilon} \left(  \frac{ e^{xt}}{t} \right)^n g_1(t) \,dt
=
\frac{ x^{n-1} e^n}{ \sqrt{ 2 \pi n}} g_1(1/x) \, \left( 1+ O(1/n) \right)
\]
uniformly on compact subsets  of the annulus $A(1/\rho;  \infty)$.
\end{lemma}
\begin{proof}
Let $x\in K \subset A(1/\rho, \infty)$.
By a change of variables, we write
\[
\frac{1}{2 \pi i}
\int_{ |t| = \epsilon} \left(  \frac{ e^{xt}}{t} \right)^n g_1(t) \,dt
=
\frac{ x^{n-1}}{ 2 \pi i}
\int_{ |t| = \epsilon |x|} \left(  \frac{ e^t}{t} \right)^n g_1(t/x) \, dt.
\]
By construction, $g(t/x)$ is analytic on  a disk of radius greater than 1.  So the contour
in the last integral can be deformed to the unit circle $|x|=1$ without changing its value.
Finally, by an application of the saddle point method we find that
\[
\frac{ x^{n-1}}{ 2 \pi i}
\int_{ |t| = 1} \left(  \frac{ e^t}{t} \right)^n g_1(t/x) \, dt
=
\frac{ x^{n-1} e^n} { \sqrt{ 2 \pi n}} g_1(1/x) \left( 1+ O(\tfrac{1}{n}) \right).
\]
\end{proof}

\begin{figure}[h]
\includegraphics[height=3.0in, width=3.0in]{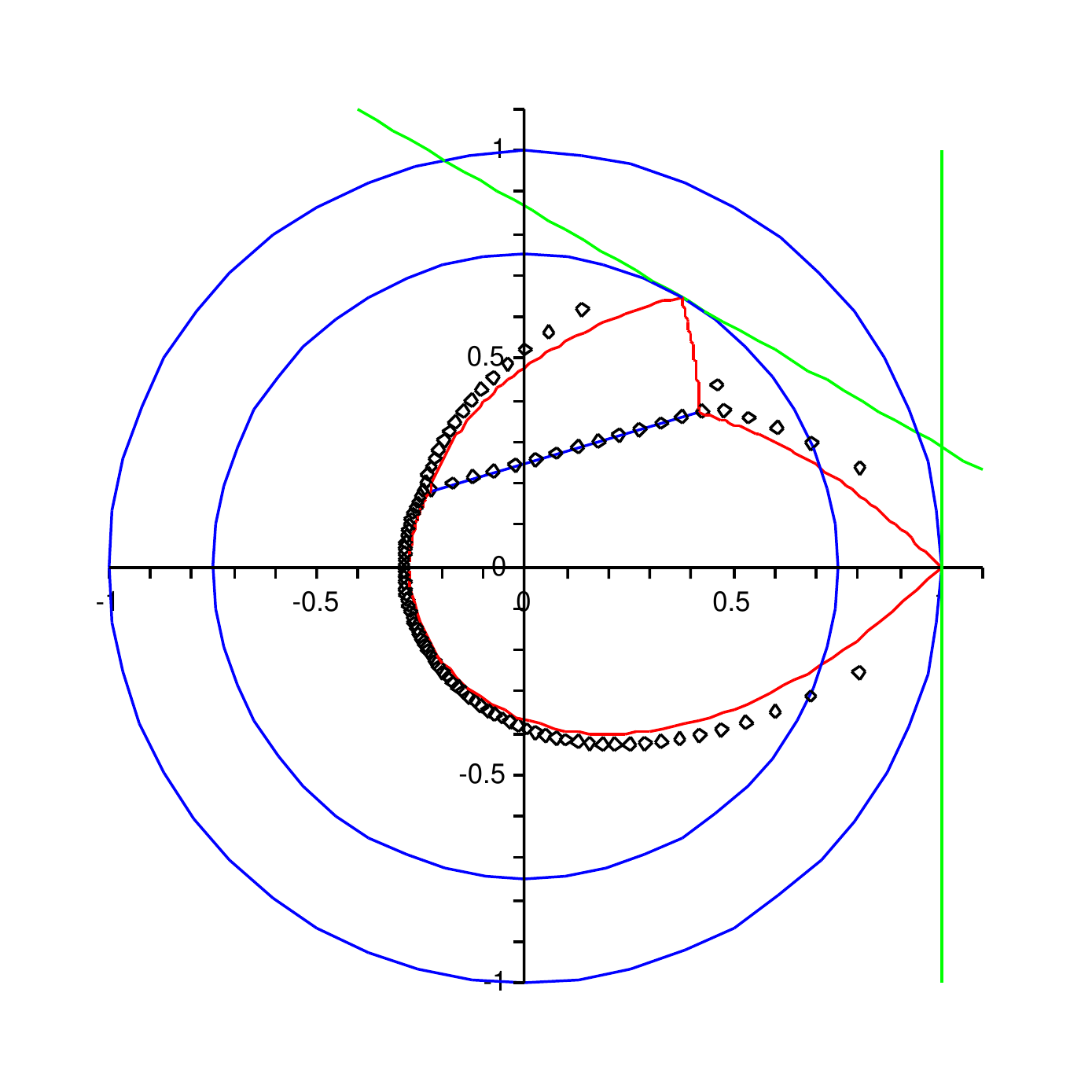} 
\caption{Generic Plot of Polynomial Zeros and Zero Attractor  When $g$ Has Two Roots; Tangent Lines and Circles Displayed }
%\label{fig:zero_attractor_1}
\end{figure}

To state the next two lemmas, we need to introduce  special polynomials $I_n(z)$ 
in $z^{-1}$ and 
$J(a; z)$ in $z$.

The polynomial $I_n(z)$ comes from expanding the derivative
of $D_{z}^{m-1}(z^{-n}e^{nz})$. Consider
\begin{eqnarray*}
D_{z}^{m-1}(z^{-n}e^{nz})
&=&\sum_{p=0}^{m-1}\binom{m-1}{p}
(D_{z}^{p}z^{-n})(D_{z}^{m-1-p}e^{nz})
\\
&=&
\sum_{p=0}^{m-1}\binom{m-1}{p}(-n)(-n-1)\cdots
(-n-p+1)z^{-n-p}(n^{m-1-p}e^{nz})
\\
&=&
z^{-n}e^{nz}n^{m-1}\sum_{p=0}^{m-1}\binom{m-1}{p}(-n)(-n-1)\cdots
(-n-p+1)(nz)^{-p}
\\
&=&
z^{-n}e^{nz}n^{m-1}\sum_{p=0}^{m-1}(-1)^{p}  p! \binom{m-1}{p}
\binom{n+p-1}{p}  \, (nz)^{-p}
\\
&=&
z^{-n}e^{nz}n^{m-1}I_{m-1}(nz),
\end{eqnarray*}
where $I_{m-1}(z)$ is given in

\begin{definition}\label{def:I}
\begin{equation}
I_{m-1}(z)=\sum_{p=0}^{m-1}(-1)^{p}  p!  \binom{m-1}{p}\binom{n+p-1}{p}
\,  z^{-p}. 
\end{equation}
\end{definition}

For $a\in Z(g)$, we define   $J(a;z)$ which are  also polynomials in $z$.
 We write out  the singular part 
$s_{a}(t)$ of the function $\frac{1}{tg(t)}$ at its nonzero pole $a$
 by
\begin{equation}\label{eq:s_a}
s_{a}(t):=\sum_{m=1}^{ \beta_a }\frac{b_{a,m}}{(t- a)^{m}},
\end{equation}
where  $\beta_a$ is the order of $a$ as a zero of $g(t)$
so
$b_{a, \beta_a} \neq 0$.

\begin{definition}\label{def:J_a}
For $a\in Z(g)$, let $J(a;z)$ be the polynomial in $z$ given as
\begin{equation}
J(a;z)=\sum_{m=1}^{ \beta_a}
\frac{b_{a,m}}{(m-1)!} \,z^{m-1}  I_{m-1}(a z).
\label{def_J}
\end{equation}
\end{definition}

\begin{lemma} \label{lemma:4.2}
Let $a \in Z(g)$ and let $x \in K$, a compact subset of the half-plane $H_a$, $\Re(ax)<1$.
Then 
\[
\frac{1}{2 \pi i}
\int_{ |t| = \epsilon} \left(  \frac{ e^{xt}}{t} \right)^n s_a(t) \,dt
=
- a^{-n} e^{n ax} J(a; nx) +
\frac{ e^n x^{n-1}}{ \sqrt{ 2 \pi n}} s_a(1/x) \left( 1+ O( n^{1-3 \alpha}) \right)
\]
where $s_a(t)$ is the singular part of $1/(tg(t))$ at the zero $a$ of $g(t)$.
\end{lemma}
\begin{proof}
We first write out the integral in terms of the singular part $s_a(t)$
\[
\frac{1}{2 \pi i}
\int_{ |t| = \epsilon} \left(  \frac{ e^{xt}}{t} \right)^n s_a(t) \,dt
=
- \sum_{m=1}^{\beta_a} \frac{ b_{a,m}}{ (m-1)!} D^{m-1}_a
\left( a^{-n} S_{n-1}(n ax) \right)
\]
where the coefficients $b_{a,m}$ are given in equation (\ref{eq:s_a}).
We now study the
asymptotics of the typical term $D_{a}^{m-1}(a^{-n}S_{n-1}(n ax))$.

 We may use  the generalized half-plane Szeg\"o asymptotics with
 $\frac{1}{3}<\alpha <\frac{1}{2}$ because of the restriction that
 $a \in Z(g)$ with $|a|\leq r_\ell$ to obtain
\begin{eqnarray*}
\lefteqn
{
D_{a}^{m-1}(a^{-n}S_{n-1}(n ax))
=
x^{n+m-1}D_{ax}^{m-1}((a x)^{-n}S_{n-1}(n ax ))
} \\
&& = \quad
x^{n+m-1}
\left\{     
\left. 
D_{z}^{m-1}(z^{-n}e^{nz})\right|_{z=ax}
 -
\frac{(m-1)!}{\sqrt{2\pi n}}\frac{e^{n}}{(1-a x)^{m}}
\left(1+O(n^{1-3\alpha }) 
\right)
\right\}   .
\end{eqnarray*}
Combining these estimates we obtain
\begin{eqnarray}
\lefteqn
{
D_{a}^{m-1}( a^{-n}S_{n-1}(n ax))
=
x^{n+m-1}
\{(ax)^{-n}e^{n ax}n^{m-1}I_{m-1}(n ax)
}  \nonumber \\
&&\qquad \quad
- \,
\frac{(m-1)!}{\sqrt{2\pi n}}\frac{e^{n}}{(1-ax)^{m}}
\left(1+O(n^{1-3\alpha }) \right)  \}
\nonumber \\
&& 
=a^{-n}e^{n ax}(nx)^{m-1}I_{m-1}(n ax)
-
\frac{(m-1)!}{\sqrt{2\pi n}}
\frac{e^{n}x^{n+m-1}}
{(1- a x)^{m}}
\left(1+O(n^{1-3\alpha }) \right).
\label{dsn-1}
\end{eqnarray}
Hence after summation we obtain
\begin{eqnarray}
\frac{1}{2\pi i}\oint_{\left| t\right| =\epsilon }
\left(\frac{e^{xt}}{t}\right)^{n}s_{a}(t)\,dt
&=&
-
\sum_{m=1}^{\beta _{a}}\frac{b_{a,m}}{(m-1)!}
D_{ a}^{m-1}( a^{-n}S_{n-1}(n ax)) 
\nonumber \\
&=&
-
a^{-n}e^{n ax}J(a;nx)+\frac{e^{n}x^{n-1}}{\sqrt{2\pi n}}
s_{a}(\frac{1}{x}) \, \left(1+O(n^{1-3\alpha }) \right). \label{eq:jzij}
\end{eqnarray}
\end{proof}

\begin{corollary} \label{cor:4.3}
For $a \in Z(g)$, $|a| \leq r_\ell$, we have
\[
\frac{1}{2 \pi i}
\int_{ |t| = \epsilon} \left(  \frac{ e^{xt}}{t} \right)^n s_a(t) \,dt
=
- a^{-n} e^{n ax} J(a; nx) +
\frac{ e^n x^{n-1}}{ \sqrt{ 2 \pi n}} s_a(1/x) \left( 1+ O( n^{1-3 \alpha}) \right)
\]
uniformly on the compact subsets of $R_\ell$,
where $s_a(t)$ is the singular part of $1/(tg(t))$ at the zero $a$ of $g(t)$.
\end{corollary}

\begin{lemma} \label{lemma:4.4}
Let $a \in Z(g)$ and let $x \in K$, where $K$ compact subset of the disk-complement
$A(1/|a|;\infty)$.
Then 
\[
\frac{1}{2 \pi i}
\int_{ |t| = \epsilon} \left(  \frac{ e^{xt}}{t} \right)^n s_a(t) \,dt
=
\frac{ e^n x^{n-1}}{ \sqrt{ 2 \pi n}} s_a(1/x) \left( 1+ O( n^{1-3 \alpha}) \right).
\]
\end{lemma}
\begin{proof}
We will use the  disk-complement generalized Szeg\"o asymptotics.
For  $z$  in the annulus $A(1+c, \infty)$, for any $c>0$, we have
\begin{equation*}
S_{n-1}(n z)=
- \frac{z^{n}}{2\pi i}
\oint_{\left| \zeta \right| =1}
\frac{
e^{n(\zeta -\ln \zeta )}}{\zeta -z} \, d\zeta
\end{equation*}
By Dividing $z^{n}$ and taking derivatives up to order $m-1$, we get 
\begin{eqnarray*}
D_{z}^{m-1}(z^{-n}S_{n-1}(nz))
&=&- \frac{ (m-1)!}{2\pi i}
\oint_{\left| \zeta \right| =1}
\frac{e^{n(\zeta -\ln \zeta )}}{(\zeta -z)^{m}} \, d\zeta
\\
&=&
-\frac{ (m-1)!}{\sqrt{2\pi n}}\frac{e^{n}}{(1-z)^{m}}
\left(1+O(n^{1-3\alpha})\right).
\end{eqnarray*}
In the above, replace $z$ by $ax$ to obtain
\begin{eqnarray}
D_{a }^{m-1}(a^{-n}S_{n-1}(n ax))
&=&
x^{n+m-1}D_{ax}^{m-1}(( a x)^{-n}S_{n-1}(n ax))
\nonumber \\
&=&
- \frac{(m-1)!}{\sqrt{2\pi n}}
\frac{e^{n}x^{n+m-1}}{(1- ax)^{m}}
\left(1+O(n^{1-3\alpha })\right).  \label{dsn-2}
\end{eqnarray}
By summation, we obtain the asymptotics for the original integral:
\begin{eqnarray}
\frac{1}{2\pi i}\oint_{\left| t\right| =\epsilon }
\left(\frac{e^{xt}}{t}\right)^{n}s_{a}(t)\,dt
&=&
-
\sum_{m=1}^{\beta_a }\frac{b_{a,m}}{(m-1)!}
D_{a}^{m-1}(a^{-n}S_{n-1}(n x))
\nonumber \\
&=&\frac{e^{n}x^{n-1}}{\sqrt{2\pi n}}  s_{a}(\frac{1}{x})
\left(1+O(n^{1-3\alpha }) \right).
\label{sij}
\end{eqnarray}
\end{proof}

\begin{corollary} \label{cor4.4}
For $a\in Z(g)$ with $ r_{\ell+1} \leq |a| < \rho$, we have
\[
\frac{1}{2 \pi i}
\int_{ |t| = \epsilon} \left(  \frac{ e^{xt}}{t} \right)^n s_a(t) \,dt
=
\frac{ e^n x^{n-1}}{ \sqrt{ 2 \pi n}} s_a(1/x) \left( 1+ O( n^{1-3 \alpha}) \right),
\]
uniformly on the compact subsets of $R_\ell$.
\end{corollary}

The remaining case for the above integration involving $s_a(t)$ on the disk $D(1/a;\delta)$ will
be handled in a later section.

\begin{proposition} \label{4.5}
For $x \in R_\ell$, we have
\begin{eqnarray*}
\lefteqn
{
\frac{ p_n(nx)}{ (xe)^n/ \sqrt{2 \pi n}}
=
\frac{1}{x} \, \frac{1}{ g_1(1/x)} 
-
\sqrt{2 \pi n}
\sum \left\{ J(a;nx) \frac{1}{ \phi(ax)^n} : a \in Z(g), |a| \leq r_\ell \right\} 
}\\
&&\qquad\qquad
+ 
\sum\left\{ \frac{1}{x} s_a\left( \frac{1}{x} \right) : a\in Z(g), |a| < \rho      \right\}       (1+    O(n^{1-3\alpha}))
\end{eqnarray*}
uniformly on the compact subsets of $R_\ell$,
where $\phi(x)= xe^{1-x}$ and $1/3<\alpha<1/2$.
\end{proposition}
\begin{proof}
Putting the last two corollaries into Equation (\ref{eq:pnnx-1}) and simplifying, we have
\begin{eqnarray*}
\lefteqn
{
p_{n}(nx)
=
\frac{x^{n-1}e^{n}}{\sqrt{2\pi n}}g_{1}(\frac{1}{x})\,
\left(1+O(\frac{1}{n})\right)
} 
\\
&& \qquad
-
\sum\left\{
a^{-n}e^{n ax}J(a; nx)  : a \in Z(g), |a| < \rho \right\}
\left(1+O(n^{1-3\alpha }) \right).
\end{eqnarray*}
 \end{proof}

\begin{proposition}  \label{prop:4.6}
For $x \in R_\ell$, we have
\[
\frac{ p_n(nx)}{ (xe)^n/ \sqrt{2 \pi n}}
=
\frac{1}{ g(1/x)} -
\sqrt{2 \pi n}
\sum\{ J(a;nx) \frac{1}{ \phi(ax)^n} : a \in Z(g), |a| \leq r_\ell \} + O(n^{1-3\alpha})
\]
uniformly on the compact subsets of $R_\ell$,
where $\phi(x)= xe^{1-x}$ and $1/3<\alpha<1/2$.
\end{proposition}
\begin{proof}

By the definition of $g_1(t)$  (see equation (\ref{eq:def_g_1})),  we see that 
\begin{equation}
\frac{1}{x} g_1 \left( \frac{1}{x} \right)
=
\frac{1}{g(1/x)}-
\sum
\left\{  \frac{1}{x}  s_{a}(\frac{1}{x})  :
a \in Z(g), |a| < \rho
\right\}.
\label{xg1}
\end{equation}
Insert this into the above. Since the $s_{a}(1/x)$ term cancels, 
we have uniformly for $x\in R_{\ell}$:
\begin{eqnarray}
\lefteqn
{
\frac{p_{n}(nx)}{(ex)^{n}/\sqrt{2\pi n}}
=
\frac{1}{g(1/x)}
}
\\ 
&& \qquad 
-\sqrt{2\pi n}
\sum
\left\{ (axe^{1- ax})^{-n}J(a; nx) :  a \in Z(g), |a| < \rho \right\}
+O\left(n^{1-3\alpha }\right).
\nonumber
\label{eq:pnnx-2}
\end{eqnarray}
\end{proof}

\begin{lemma} \label{lemma:4.7}
\label{lemma:jzijnx} If $a \in Z(g)$ with $|a| < \rho$ and $ x\in R_{\ell}$,
then
\begin{equation*}
J(a; nx)=\frac{  b_{a,\beta_a}
}
{( \beta_a-1)!}
(nx)^{ \beta_a -1}
\left(\frac{a x-1}{a x}\right)^{\beta_a -1}
(1+o(1)).
\end{equation*}
\end{lemma}

\begin{proof}  
 Recall that
\begin{equation*}
J(a; nx)=\sum_{m=1}^{\beta_a}
\frac{b_{a,m}I_{m-1}(n a x)}{(m-1)!}, 
\quad
I_{m-1}(n  a x)=\sum_{p=0}^{m-1}(-1)^{p}
\binom{m-1}{p}\binom{n+p-1}{p} p! \, (n a x)^{-p}.
\end{equation*}
It is easy to see that
\begin{equation*}
\binom{n+p-1}{p}(n a x)^{-p}
=
\frac{( a x)^{-p}}{p!} \left(1+o(1) \right),
\end{equation*}
that is, as $n\rightarrow \infty $
\begin{equation*}
I_{m-1}(n a x)\rightarrow 
\sum_{p=0}^{m-1}(-1)^{p}\binom{m-1}{p}
(a x)^{-p}= \left(\frac{ a x-1}{a x} \right)^{m-1}.
\end{equation*}
Hence
\begin{equation}
J(a;nx)
=
\frac{  b_{a,\beta_a}
}
{(\beta_a-1)!}
(nx)^{\beta_a -1}
\left(\frac{a x-1}{ a x}\right)^{\beta_a -1}
(1+o(1)).
\label{esti-j}
\end{equation}
Since the coefficient $b_{a,\beta_a}$ in the definition of the
singular part $s_{a}(t)$ is nonzero, 
we find for fixed $x$ that  the precise order of $J(a;nx)$ as a polynomial in $n$ is $n^{ \beta_a-1}$.
\end{proof}

We note the following
\begin{corollary} \label{cor:4.8}
$\displaystyle
\lim_{n\to \infty} \frac{ p_n(nx)}{ (xe)^n/ \sqrt{2 \pi n}} = \frac{1}{g(1/x)}$, $x \in R_\ell
$
provided $| \phi(ax) | >1$ for all $a \in Z(g)$ with $|a| \leq r_\ell$.
\end{corollary}

\begin{corollary} \label{cor:4.9}
$\displaystyle
\lim_{n\to \infty} 
\frac{1}{n} \ln \left|  \frac{ p_n(nx)}{ (xe)^n/ \sqrt{2 \pi n}} \right| = - \ln | \phi(a_0 x)|$, for all   
$x \in R_\ell$
satisfying  $| \phi(a_0 x) | >1$ for  $a_0\in Z(g)$ and $|\phi_0(ax)| \leq | \phi(ax)| $ for all $a \neq a_0$
such that $|a| \leq r_\ell$.
\end{corollary}

By construction, $R_\ell \subset R_\rho$ for all $\ell$ chosen so $r_\ell < \rho$. Consequently,
we have the asymptotics:

\begin{theorem} \label{thm:4.10}
On $R_\rho$, we have the following uniform asymptotics
\begin{eqnarray*}
\lefteqn
{
\frac{p_{n}(nx)}{(ex)^{n}/\sqrt{2\pi n}}
=
\frac{1}{g(1/x)}  \, \left(1+ O(1/n) \right)
}
\\ 
&& \qquad 
-\sqrt{2\pi n}
\sum
\left\{ (axe^{1- ax})^{-n}J(a; nx) :  a \in Z(g), |a| < \rho \right\}
+O\left(n^{1-3\alpha }\right),
\end{eqnarray*}
where $1/2 < \alpha < 1/3$.
\end{theorem}

It remains to develop the asymptotics in the disks $D(1/a;\delta_a)$ and well as determining domination 
among $a\in Z(g)$ of $| \phi(ax)|$.

\section{Dominant Zeros and Szeg\"o Curves}  % section 5

Let $\phi(x)=x e^{1-x}$ which is an entire function that
is conformal on the unit disk.
The standard Szeg\"o curve $\mathcal S$ 
is the portion of the level
curve $| \phi(x) | =1$ that lies inside the closed unit disk or
equivalently, inside the closed left-hand plane $\Re(x) \leq 1$.
$\mathcal S$ is a closed simple closed convex curve; in fact,
it has the form $t= \pm \sqrt{ e^{2(x-s)}-s^2}$ where
$x=s+it$ and $s \in [ -W(e^{-1}),1]$ and $W$ is the principal
branch of the Lambert $W$-function.

\begin{definition}
Let $a$ be a nonzero complex number.
We call any curve of the form $\frac{1}{a} {\mathcal S}$ a Szeg\"o curve.
\end{definition}

\begin{figure}[h]
\includegraphics[height=3.0in, width=3.0in]{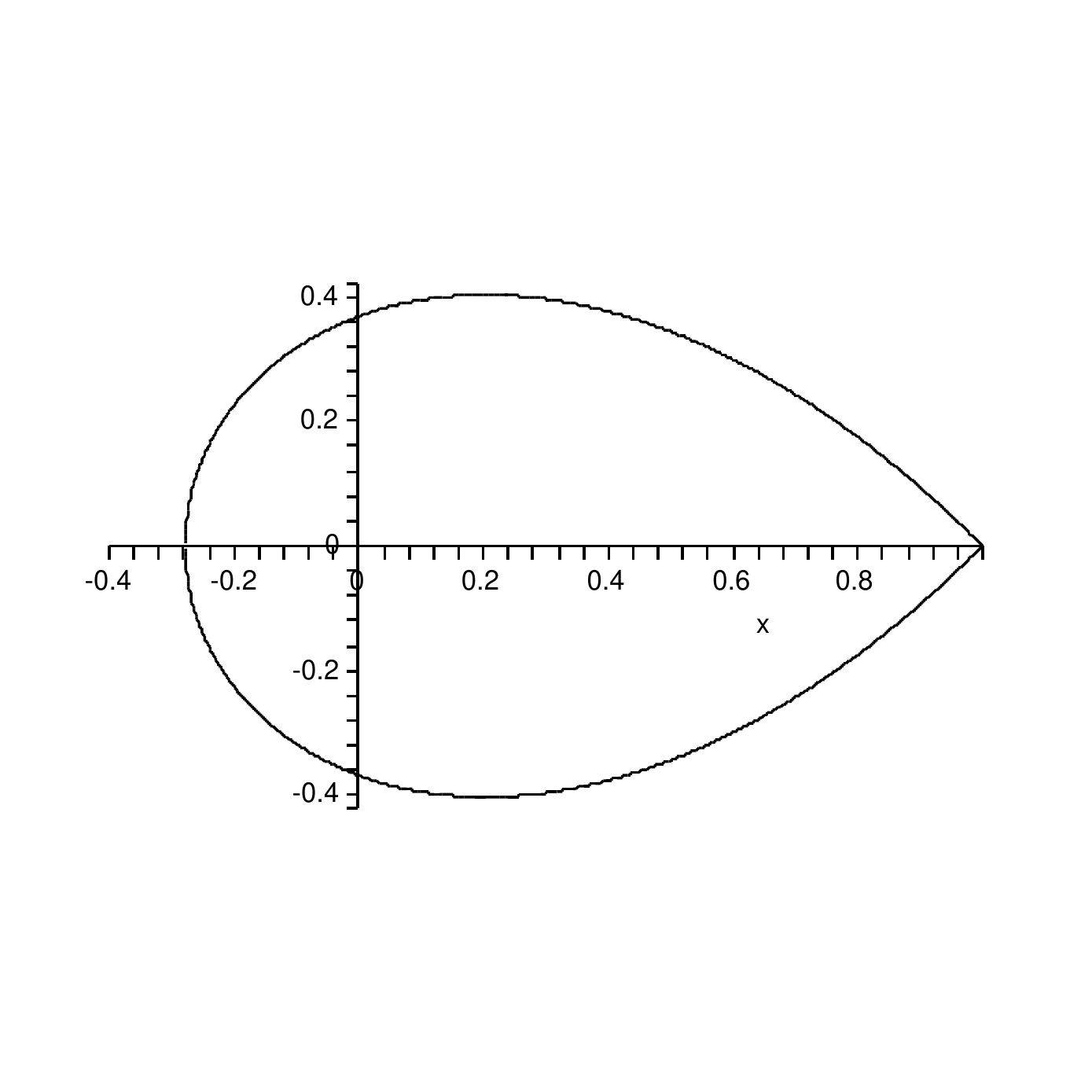} 
\caption{Szeg\"o  Curve: $ | ze^{1-z}|=1$ and $|z| \leq 1$ }
\label{fig:szego_curve}
\end{figure}

\begin{lemma} \label{lemma:5.1}
Let $a, b$ be non-zero distinct complex numbers. Then 
the intersection $\frac{1}{a} \mathcal{S} \cap \frac{1}{b} \mathcal{S}$
has at most 2 points.
\end{lemma}
\begin{proof}
The intersection of the two curves  $\frac{1}{a} \mathcal{S} \cap \frac{1}{b} \mathcal{S}$
must satisfy $| \phi(ax)| = | \phi(bx)|$. This modulus condition  determines a line; so, 
the intersection of the two curves lie on the $\frac{1}{a} \mathcal{S}$ and a line. Since
$\frac{1}{a} \mathcal{S}$ is convex, the intersection contains at most two points.

Write $x = s+it$ and $b-a= \alpha+ i \beta$. Then $ | \phi(ax)| = | \phi(bx)|$ reduces to the line:
\begin{eqnarray*}
| a e^{-ax} | =  | b e^{ -bx}|,\quad
| e^{(b-a) x} | = |b/a|,\quad
| e^{ \Re[ (b-a)x]}| = | b/a|,\\
\Re[ (b-a)x] = \ln |b/a|,\quad
\alpha s - \beta t = \ln |b/a|.
\end{eqnarray*}
\end{proof}

\begin{definition}
Call $a \in Z(g)$ a dominant zero if either $|a| = r_0$ or
$1/a$ does not lie in the interior of any of Szeg\"o curves $\frac{1}{b} {\mathcal S}$
where $b$ is a minimal modulus zero, $|b|=r_0$.
If $1/a$ lies on one of these Szeg\"o curves, call $a$ an improper dominant zero;
otherwise, $a$ is proper.
\end{definition}

Note: if $a$ is an improper dominant zero, then 
there exists a minimal modulus zero $b$ so
the intersection 
$\frac{1}{a} {\mathcal S} \cap \frac{1}{b} {\mathcal S}$ consists of
a single point, namely $1/a$.
In \cite{bleher},  the behavior of the zeros of 
the Taylor polynomials of $\sum_{j=1}^m c_j e^{\lambda_jz}$
is determined if the numbers $\lambda_j$ satisfy a condition in terms of
a convex $m$-gon.  For the Appell polynomials, the geometric condition
is more subtle. 

Let $W$ denote the principal value of the Lambert $W$-function. Since the radius of
the largest circle centered 0 that lies in the interior of the standard Szeg\"o curve
$\mathcal S$ is $W(e^{-1})$, we have the following

\begin{lemma}
If $a'$ is a zero of $g$ such that $| a'| > r_0/ W(e^{-1})$, then $a'$ must be a non-dominant
zero. Hence, there are at most finitely many dominant zeros.
\end{lemma}

\begin{definition}
If $a$ and $b$ are two dominant zeros such that $\frac{1}{b} {\mathcal S} \cap \frac{1}{a} {\mathcal S}$
consists of two points, then the line $| \phi(ax)|= |\phi(bx)|$ determines two half-planes. Let
$E_{a^+,b}$ denote the closed half-plane that contains $1/a$. 
\end{definition}

Let 
$a_1,a_2, \dots, a_n$ be the dominant zeros of $g$.
Recall that if $f(x)$ is any analytic function on a domain $D$
then its modulus $| f(x)|$ is a subharmonic  function on $D$.
We need two basic properties of subharmonic functions:
 they satisfy the maximum modulus principle; and
 the maximum of two subharmonic functions is still
subharmonic.

\begin{definition}\label{def:phi}
Set
$ \Phi  (x) = \max( | \phi(a_1x)|^{-1}, \cdots, | \phi(a_n x)|^{-1})$
so $\Phi (x)$ is subharmonic on the punctured complex plane.
\end{definition}

We work with  $| \phi(a x)|^{-1}$ rather than $ | \phi(a x)|$ so we may
apply the Maximum Principle since $ |\phi(a x)|$ vanishes at $x=0$.

We observe that the  level curve $\Phi(x)=1$
 divides the complex plane
 into finitely many connected components just as the original  curve $| \phi(x)|=1$
 divides the complex plane into three connected components.
 We make the:

\begin{definition}\label{def:D0}
Consider the connected components of the complement of the level curve $\Phi(x)=1$.
Let $D_0$ be the closure of the connected component that contains $0$.
Note that if $x\in D_0$ and is $x\neq 0$, then $  \Phi(x)\geq 1$ with
strict inequality when $x$ lies in the interior of $D_0$.
\end{definition}

We assume that  $\rho>0$ is chosen so large that  the closed disk $\overline{D}(0;1/\rho)$
is a subset of the interior of $D_0$.
 Furthermore,
the singularity in $\log | \phi(a x)| - \log | \phi(b x)|$
always cancels so this difference is always a harmonic function
on $\mathbb C$.

Let $a'$ be a non-dominant zero of $g$ such that $|a'| < \rho$. 
Then we know that
 $\frac{1}{a'} {\mathcal S}$ is a subset of the
interior of $D_0$ and that the disk $D(0;1/\rho)$ lies 
inside of $\frac{1}{a'} {\mathcal S}$.

The following Proposition follows easily from the definition of
$D_0$ and dominant zeros:

\begin{proposition} \label{prop:5.2}
Uniformly on the compact subsets of  $R_\rho \setminus D_0$,  we have
\[
\lim_{n\to\infty} \frac{ p_n(nx)}{ (ex)^n/\sqrt{2\pi n}} = \frac{1}{g(1/x)}
\]
\end{proposition}

To understand the asymptotics inside $D_0$ requires subharmonic
function theory.

By construction, this means that $|\phi( a' x)^{-1} | <1$
for $x \in \partial D_0$, while  $\Phi(x)=1$ for
$x \in \partial D_0$.  Before we can apply the Maximum
Modulus Principle for subharmonic functions, we need to
deal with the common singularity at 0. However, this is easily
dealt  with by multiplying both $\Phi(x)$ and $|\phi^{-1}(x)|$
by $|x|$ which shows that we can remove this singularity
relative to the inequality.

In particular, $\Phi(x)$ is strictly larger than $| \phi^{-1}(x)|$
for all $x \in D_0$, with $x \neq 0$. So there exists a positive
constant, say $\alpha$, so that
\[
 \Phi(x)  >  | \phi(a'x)|^{-1} + \alpha, \quad x \in D_0 \setminus D(0;1/\rho).
\]

We state the above discussion  formally:

\begin{proposition} \label{prop:5.3}           %\label{prop:order_non_dom}
Let $a'$ be a non-dominant zero of $g$.
On the domain $D_0 \setminus D(0;1/\rho)$, we have the order
estimates
\[
| \phi(a'x)|^{-n} 
= 
O(  ( \Phi(x) - \alpha)^n)
=o (  \Phi (x)^{n}). 
\]
\end{proposition}

\begin{definition}\label{def:D_a}
Given a proper dominant zero $a$ of $g$, let
\[
D_a 
= \{ x \in D_0 : | \phi(ax)|  \leq | \phi(bx)|, 
\textrm{  for all dominant zeros  } b \neq a,  | \phi(ax)| \leq 1 \}.
\]
\end{definition}
Note that the definition of $D_a$ is independent of the choice of $\rho$
sufficiently large. Further, $D_a$ has the alternate description 
 in terms of the half-planes $E_{a^+,b}$:
\[
\frac{1}{a} {\rm Interior}({\mathcal S} )
\cap \, \bigcap \{ E_{a^+,b} : b \textrm{  is a proper dominant zero of  } g, b \neq a \}.
\]

It is easy to state formally the basic properties of $D_a$:

\begin{lemma} \label{lemma:5.4}
Let $a$ be a proper dominant zero of $g$. Then $D_a$ is a non-empty
compact connected subset of $D_0$.
\end{lemma}
\begin{proof}
Since $\frac{1}{a}{\mathcal S}$  is a convex curve and the intersection
of half-planes is connected, the set $D_a$ must be a connected
convex set.
\end{proof}

We now restate Proposition \ref{prop:5.3} relative to a domain $D_a$:

\begin{proposition}  \label{prop:5.5}                 %         \label{prop:non_dom_order}
Let $a'$ be a non-dominant zero of $g$ with $|a| < \rho$. Then $a'$ must lie
in a domain $D_a$ for some proper dominant zero $a$ of $g$. For $x\in D_a \setminus D(0;1/\rho)$,
the following holds uniformly
\[
| \phi(a'x)|^{-n} 
= 
O(  | \phi(ax) - \alpha|^{-n})
=o (  |\phi (ax)|^{-n}  ). 
\]
\end{proposition}

Finally, we now have a refinement of Theorem \ref{thm:4.10} as

\begin{theorem}\label{thm:5.6}
Let $\rho$ be chosen greater than $1/|a|$ where $a$ is any proper dominant zero of $g$. Then
on $R_\rho$, we have the following uniform asymptotics
\begin{eqnarray*}
\lefteqn
{
\frac{p_{n}(nx)}{(ex)^{n}/\sqrt{2\pi n}}
=
\frac{1}{g(1/x)} \, \left( 1+O(1/n) \right)
}
\\ 
&& \qquad 
-\sqrt{2\pi n}
\sum
\left\{  \phi(ax)^{-n}J(a; nx) :  a \in Z(g)  \textrm{ and dominant  } \right\}
+O\left(n^{1-3\alpha }\right) + o( \Phi(x)),
\end{eqnarray*}
where $1/3<\alpha<1/2$ and $\Phi(x) = \max\{ | \phi(ax) |^{-1} : a\in Z(g)  \textrm{ and dominant  } \}$.
\end{theorem}

\section{Asymptotics Inside the Disk $D(1/a'; \delta_{a'})$ When $a'$ is a Non-Dominant Zero}   % section 6

We first state an easy consequence of a previous Proposition.

\begin{proposition} \label{prop:6.1}
Let $a' \in Z(g)$. Then on the disk $D(1/a'; \delta_{a'})$, the normalized polynomials have the asymptotics
\begin{eqnarray*}
\qquad \frac{ p_n(nx)}{ (xe)^n/\sqrt{2 \pi n}}
=
\frac{1}{x} g_1\left( \frac{1}{x} \right) \, \left( 1+ O\left( \frac{1}{n} \right) \right)
+
\sum
\left\{
J(a;nx) \frac{1}{ \phi(ax)^n} : a \in Z(g), |a| \leq |a'|, a \neq a'
\right\} 
\\
+
\sum
\left\{
\frac{1}{x} s_a\left( \frac{1}{x} \right) : a \in Z(g), a \neq a', |a| < \rho
\right\}
+
\sigma_{a'}(x), \qquad\qquad&&
\end{eqnarray*}
where 
\[
\sigma_{a'}(x) =
 \sum_{m=1}^{ \beta_{a'}} \frac{ b_{a',m}}{ (m-1)!} D^{m-1}_{a'}\left(  (a')^{-n} S_{n-1}(n a'x) \right).
\]
\end{proposition}

\begin{proposition} \label{prop:6.2}
% \label{pnnx-d} 
Let $a'$ be a non-dominant zero of $g$ with $|a'| < \rho$.
Then there exists a choice of $\delta_{a'}>0$ such that 
\[
\sigma_{a'}(x)=O\left(e^{  6n \delta_{a'} \rho } \right)
\]
where $\rho$ is the cut-off bound for the zeros of $g$.
\end{proposition}

\begin{proof}
To estimate $\sigma_a$, we make use of the elementary estimate:
If $f(z)$ is analytic function of $z$, then for any 
$\epsilon >0$, we have
\begin{equation*}
\left| D_{z}^{j-1}f(z)\right|
 \leq \frac{(j-1)!}{\epsilon^{j-1}}
\max_{\left| \zeta -z\right| 
=\epsilon }\left| f(\zeta )\right| .
\end{equation*}

By the definition of $\sigma_{a'}(x)$, we find
\begin{eqnarray*}
\left| \sigma_{a'} \right| 
&\leq&
 \left| \sum_{m=1}^{\beta_{a'} }
 \frac{b_{a',m}}{(m-1)!}D_{a'}^{m-1}( | a' |^{-n}S_{n-1}(n a' x))\right|
\\
&\leq&
 \sum_{m=1}^{\beta _{a'}}
 \frac{\left| b_{a',m}\right| }{\delta_{a'}^{m-1}}  \,
\,\max_{\left| \zeta -a' \right| =\delta_{a'} }
\left| \zeta^{-n}S_{n-1}( n \zeta  x)\right| 
\\
&\leq&
 K_{ \delta_{a'} }\,
  \max_{\left| \zeta -a' \right| =\delta_{a'} }
  (\left| \zeta\right|^{-n}
  S_{n-1}(\left| \zeta x\right| n)   )
\end{eqnarray*}
where $K_{\delta_{a'}}>0$ is a constant that depends on the zero $a'$ and the
radius $\delta_{a'}$.

To go further we observe for $x\in D(\frac{1}{a},\delta_{a'} )$
and $| \zeta - a'| = \delta_{a'}$:
\begin{eqnarray*}
\left| \zeta x\right| \leq (\left| a' \right| +\delta_{a'} )\left| x\right|
&\leq&
 \left| a' \right| \left| x\right| +\left| x\right| \delta_{a'}
\\
&\leq&
 1+\left| a' \right| \delta_{a'} +\left| x\right| \delta_{a'} =1+\delta_{a'}
(\left| a' \right| +\left| x\right| ).
\end{eqnarray*}
Since  $\left| a' \right| <   \rho$ by assumption,
$\left| \zeta x\right| \leq 1+ 2 \rho \delta_{a'}$.
But  $\left| \zeta \right| \geq \left| a' \right| - \delta_{a'}$ and
$\left| x\right| \geq \frac{1}{\left| a' \right| }-\delta_{a'}$, so we get
\begin{equation*}
\left| \zeta x\right| \geq \left(\left| a' \right| - \delta_{a'} \right)
\,
\left(\frac{1}{\left| a' \right| }-\delta_{a'} \right)
\geq 1-
\delta_{a'} 
\left(
\frac{1}{\left| a' \right| }
+
\left| a' \right| 
\right)
   \geq 1-2\delta_{a'}  \rho.
\end{equation*}
Collecting these two inequalities, we get 
\begin{equation*}
1- 2\delta_{a'} \rho \leq \left| \zeta x\right| \leq 
1+   2\delta_{a'}  \rho.
\end{equation*}
Now use that $\left| S_{n-1}(nt)\right| \leq e^{nt}$:
\begin{eqnarray*}
\max_{\left| \zeta - a' \right| =\delta_{a'} }
\left| e\zeta x\right|^{-n}S_{n-1}(\left| \zeta x\right| n)
&\leq& 
e^{-n}\left| 1-  2\delta_{a'} \rho \right|^{-n}   e^{n(1+ 2\delta_{a'}  \rho)}
\\
&=&
\left| 1- 2\delta_{a'}  \rho\right|^{-n}
e^{  2 n \delta_{a'}  \rho}
\end{eqnarray*}
For $0\leq x\leq {1}/{2}$, ${1}/{(1-x)}\leq e^{2x}$;
if we choose $\delta_{a'} $ such that
$
2\delta_{a'}  \rho \leq {1}/{2},
$
then we have $\left| 1- 2\delta_{a'}  \rho\right|^{-n}
\leq 
e^{ 4\delta_{a'} \rho }$. With this choice of $\delta$,
we obtain the desired bound
\begin{equation*}
\max_{\left| \zeta - a' \right| =\delta_{a'} }
(\left| e\zeta x\right|^{-n}S_{n-1}(\left| \zeta x\right| n))
\leq e^{ 4\delta_{a'} \rho}
e^{ 2n \delta_{a'}  \rho }=e^{6n \delta_{a'}  \rho}
\end{equation*}
\end{proof}
\section{Zero Attractor and the Density of the Zeros}   % section 7

In our paper \cite{g-b}, we determined the limit points of the zeros
of the Euler polynomials by means of the asymptotics and the zero
density.  Here, we separate out first the question of find the support
of the zero density measure, which is, of course, the zero attractor.
Then we determine the zero density by applying our general result
in the appendix.

\begin{proposition} \label{prop:7.1}
Let $f_n(x) = \sqrt{2\pi n} \,p_n(nx)/ (xe)^n$. Then the following limits hold uniformly
on compact subsets
of  the indicated domains:
\begin{enumerate}
\item
On the domain $A(1/r_0; \infty)$,
$\displaystyle
\lim_{n\to\infty} \frac{1}{n}  \ln[ f_n(x)]
= 0
$.
\item
On the domain $R_\rho $,
$\displaystyle
\lim_{n\to\infty} \frac{1}{n}  \ln[ f_n(x)]
= 0
$.
\item
On the domain $D_a \cap A(1/\rho;\infty)$ where $a$ is any dominant zero of $g$,
$\displaystyle
\lim_{n\to\infty} \frac{1}{n}  \ln[ f_n(x)]
=  - \ln \phi(ax)$.
\end{enumerate}
\end{proposition}
\begin{proof}
We use the asymptotic expansions for $p_n(nx)$ developed in the previous sections.
For $|x|>1/r_0$, we noted already that the indicated limit must be 0.

We observe that if $a'$ is a nondominant zero of $g$ with $|a'| <\rho$,  then for $\delta_{a'}>0$ sufficiently
small, the disk $D(1/a'; \delta_{a'})$ will lie in the domain $D_a$ for some dominant zero
$a$, then on $D(1/a'; \delta_{a'})$,
$\displaystyle
\lim_{n\to\infty} \frac{1}{n}  \ln[ f_n(x)]
=  - \ln \phi(ax)$.
\end{proof}

To describe the zero attractor requires a closer examination of the boundary of
each domain $D_a$ where $a$ is a proper dominant zero.

The boundary $\partial D_a$ where $a$ is a proper dominant zero of $g$ has several
natural families: $\partial D_a \cap \partial D_0$ which is an ``outer boundary" and
a polygonal curve consisting of
the  line segments $\partial D_a \cap \partial D_b$ where $b$ is another
dominant zero of $g$. Note that $\partial D_a \cap \partial D_b$ is a subset of $D_0$.
It will be useful to subdivide $\partial D_a \cap  \partial D_0$ into two connected components
denoted by $\partial D_a^{\pm}$  that come
from deleting $\{1/a \}$ from $[\partial D_a \cap \partial D_0]$.

\begin{lemma}
The zero attractor of the Appell polynomials  $\{ p_n(nx)\}$ must lie inside the compact set
\[
\bigcup \left\{
\partial D_a : a \textrm{  is a proper dominant zero of  }  g 
\right\}.
\]
\end{lemma}
\begin{proof}
First, we let $x^*$ let in the infinite exterior of $D_0$. Recall that 
$\lim_{n\to\infty} \sqrt{2\pi n} \, p_n(nx)/ (xe)^n = 1/g(1/x)$ uniformly on compact subsets. 
If $x_{n_k}$ is a zero of $p_{n_k}(n_k x)$ and $x_{n_k} \to x^*$, then appealing to this limit
we find that the limit must be 0 while the right-hand side is $1/g(1/x^*) \neq 0$. 
Secondly, suppose $x^*$ lies in the interior of $D_0$ but not on any boundary set $\partial D_a$,
where $a$ is a dominant zero. By construction, $x$ will lie in the interior of one of the domains 
$D_b$, where $b$ is a dominant zero. 
Then $\lim_{n\to\infty} | \sqrt{2\pi n} \, p_n(nx)/ (xe)^n |^{1/n} = | \phi(bx)|$
uniformly on compacta in the interior of $D_a$.
By the same reasoning
as before, $x^*$ cannot be a limit of zeros.
\end{proof}

The following Theorem is an immediate consequence of the above lemma
together with the result of Sokal in section \ref{thm:sokal} of the Appendix.

\begin{theorem} \label{thm:7.2}
The zero attractor of the Appell polynomials  $\{ p_n(nx)\}$ is given by
\[
\bigcup \left\{
\partial D_a : a \textrm{  is a proper dominant zero of  }  g 
\right\}.
\]
\end{theorem}
\begin{proof}
Let $a$ be any proper ominant zero of $g$ and
let $x^*\in \partial D_a^{\pm}$. Let $\epsilon>0$ be given.
Then we find that
\begin{eqnarray*}
\lim_{n\to\infty} \ln
\left|  \frac{ p_n(nx)}{ (xe)^n / \sqrt{2\pi n}}  \right|
=
\left\{
\begin{array}{rl}
0,  & x\in D(x^*;\epsilon) \setminus D_0,\\
- \ln| \phi(ax)|, & x \in D(x^*;\epsilon) \cap {\rm Int}(D_0)
\end{array}\right.
\end{eqnarray*}
holds uniformly on compact subsets.
Next suppose that $x^*$ is nonzero and  lies on one of the line segments of the
form $\partial D_a \cap \partial D_b$ where $b$ is another proper dominant zero.
Again, we find that
\begin{eqnarray*}
\lim_{n\to\infty} \ln
\left|  \frac{ p_n(nx)}{ (xe)^n / \sqrt{2\pi n}}  \right|
=
\left\{
\begin{array}{rl}
- \ln| \phi(ax)|  ,  & x\in D(x^*;\epsilon) \cap {\rm Int}(D_a),\\
- \ln| \phi(bx)|, & x \in D(x^*;\epsilon) \cap {\rm Int}(D_b)
\end{array}\right.
\end{eqnarray*}
which also holds uniformly on compact subsets.
By Sokal's result \cite{sokal} which is described in the appendix, we conclude that $x^*$ is in $\limsup Z(p_n)$
since there can be no harmonic function $v(x)$ on the disk $D(x^*;\epsilon)$ that
satisfies the inequalities
\[
\liminf_{n\to\infty} \ln
\left|  \frac{ p_n(nx)}{ (xe)^n / \sqrt{2\pi n}}  \right| \leq v(x)
\leq
\limsup_{n\to\infty} \ln
\left|  \frac{ p_n(nx)}{ (xe)^n / \sqrt{2\pi n}}  \right|.
\]

On the disks $D(1/a'; \delta_{a'})$ where $a'$ is a non-dominant zero
which must lie inside $D_0$,
the contribution of the dominant zero dominants.

This reasoning handles all but finitely many points: $1/a$ where $a$ is
a dominant zero of $g$. However, since the zero attractor must be
a compact set and points in $D(1/a;\epsilon) \cap [ \partial D_a \cap \partial D_0]$
lie in the zero attractor, we conclude that $1/a$ also lie in the attractor.
\end{proof}

\begin{theorem} \label{thm:7.4}
Let $a$ be a proper dominant zero of $g$.
\\
(a)
The zero density measure on any proper subcurve of $\partial D_a \cap \partial D_0$ is the pull-back of
the normalized Lebesgue measure on the unit circle under the conformal map
$\phi(ax)$.
\\
(b)
Let $b$ be a proper dominant zero of $g$ so $b\neq a$. Then
 the zero density measure on any proper  line segment of $\partial D_a \cap \partial D_b$ 
 is Lebesgue measure.
\end{theorem}
\begin{proof}
For both parts, we can use the asymptotics given in Theorem \ref{thm:5.6}.

For part (a), let $f_n(x) = \sqrt{2\pi n} g(1/x) p_n( nx)/ (xe)^n$. Let $a$ be a dominant zero of $g$,
and
let $C$ be a proper subcurve of $\partial D_0 \cap \partial D_{a}^{\pm}$.
Then  there exists a neighborhood $U$ of $C$
such that
$U \subset R_\rho \cap  [( {\mathbb C} \setminus D_0 )  \cup D_a]$ so that 
the asymptotics in Theorem \ref{thm:5.6} can be written as
\[
\frac{ p_n(nx)}{ (xe)^n/ \sqrt{2\pi n}} =
\frac{1}{g(1/x)}\, \left( 1+O(1/n) \right)
- \sqrt{2 \pi n}\, \frac{ J(a;nx)}{ \phi(ax)^n} + O(n^{1-3\alpha}) + o(\Phi_{1,a}^n (x)),
\]
where $ \Phi_{1,a}(x) = \max\{ 1, | \phi(ax)| \}$.
Hence, by dividing by $g(1/x))$, we find that $f_n(x)$ has the form:
\[
f_n(x)=
1+ a_n(x) \phi(a x)^{-n} + e_n(x), \quad a_n(x) = - g(1/x) J(a;nx),
\]
where
\[
e_n(x)=
\left\{
\begin{array}{rl}
o(1),& \, x \in U \cap ({\mathbb C} \setminus D_0),\\
o( \phi(a x)^{-n}),& \, x\in  U\cap  (D_{a}  \cap R_\rho) .
\end{array} \right.
\]
Since $\phi(a x)$ is conformal in the disk $D(0;1/|a|)$, we may
apply Theorem \ref{thm:density} from the Appendix Section \ref{thm:sokal}
on the density of zeros.

Let $a$ and $b$ be two distinct proper dominant zeros of $g$ such that
$\partial D_a \cap \partial D_b$ is nonempty. On $D_a \cap D_b \cap R_\rho$,
the asymptotics in Theorem \ref{thm:5.6} can be written as
\begin{eqnarray*}
\lefteqn
{
\frac{ p_n(nx)}{ (xe)^n/\sqrt{2\pi n}}
=
\frac{1}{g(1/x)} \left( 1+O(1/n) \right) - \sqrt{2 \pi n}
\left(
J(a;nx)\frac{1}{\phi(ax)^n} + J(b;nx)\frac{1}{\phi(bx)^n} \right.
} \\
&&\qquad \left.
+ \sum \{ J(a';nx) \frac{1}{ \phi(a'x)^n} : a' \textrm{ proper dominant zero}, a' \neq a,b\,\}
\right) + O(n^{1-3\alpha}) + o( \Phi(x)^n)
 \\
&&=
\frac{1}{g(1/x)} \left( 1+O(1/n) \right) - \sqrt{2 \pi n}
\left(
J(a;nx)\frac{1}{\phi(ax)^n} + J(b;nx)\frac{1}{\phi(bx)^n}
\right) + O(n^{1-3\alpha}) + o( \Phi_{a,b}^n(x)),
\end{eqnarray*}
where $\Psi_{a,b}(x) = \max\{ 1/|\phi(ax)|, 1/ |\phi(bx)|\}$.

Let $L$ be a proper line segment of
the intersection $\partial D_a \cap \partial D_b$.
Let $U$ be a neighborhood of $L$ so both
$|\phi(ax)|<1$ and $ | \phi(bx)|<1$ for $x\in U$. 
On the intersection
 $U \cap R_\rho$, we work with a different normalization than before:
 \[
T_n(x) =  - \frac{ \phi(ax)^n}{ \sqrt{2\pi n} (xe)^n J(a;nx)} \, p_n(nx).
 \]
 Note that in this normalization the term that contains $\phi(ax)^{-n}$ becomes the
 constant 1 for $T_n(x)$.
 Of course, this new normalization has exactly the same zeros as $p_n(nx)$
 in $U$ so the zero density is unchanged.
 Then we find that
 \[
 T_n(x)
 = 1 + a_n(x) \psi(x)^n + e_n(x),
 \]
 where
 \[
 \psi(x) = \frac{ \phi(ax)}{ \phi(bx)} =\frac{a}{b}e^{ (b-a) x}, \quad a_n(x) = \frac{ J(b;nx)}{J(a;nx)},
 \]
 and
 \[
 e_n(x)= - \frac{ \phi(ax)^n}{ \sqrt{2\pi n} J(a;nx)} \left(
 O(n^{1-3\alpha}) + o(\Phi_{a,b}^n(x)) \right).
  \]
  On $U$,  we have that $ \phi(ax)^n \Phi_{a,b}^n(x) =  \max\{ 1, | \psi(x)|^n\}$; 
  while on $D_a \cap U$, $| \psi(x)|<1$ and on $D_b\cap U$, $|\psi(x)|>1$.
  This allows us to write $e_n(x)$ as
 \[
  e_n(x)=
  \left\{
  \begin{array}{rl}
  o( \psi(x)^n) ,&\, x\in D_a \cap U,\\
  \\
o(1) ,&\, x\in D_b \cap U.
  \end{array} \right.   .
 \]
 By construction,  $\phi(ax)/\phi(bx) = \frac{a}{b}e^{ (b-a) x}$ is a conformal
 map on $U\cap R_\rho$ that maps $L$ onto an arc of the unit circle. 
  By Corollary \ref{cor:constant} in the Appendix section \ref{thm:sokal}, the result follows.
\end{proof}

We close with several examples that illustrate the main constructions in the paper.

\begin{example}
{\rm
Let $g(t)$ be an entire function whose minimal modulus zero $a_1=1$ such that all its
other zeros $a$ satisfy $1/|a| < W(e^{-1}) \simeq 0.27846$. Then the zero attractor for the
associated Appell polynomials coincide with the classical Szeg\"o curve in
Figure \ref{fig:szego_curve}.
}
\end{example}

\begin{figure}[h]
\includegraphics[height=3.0in, width=3.0in]{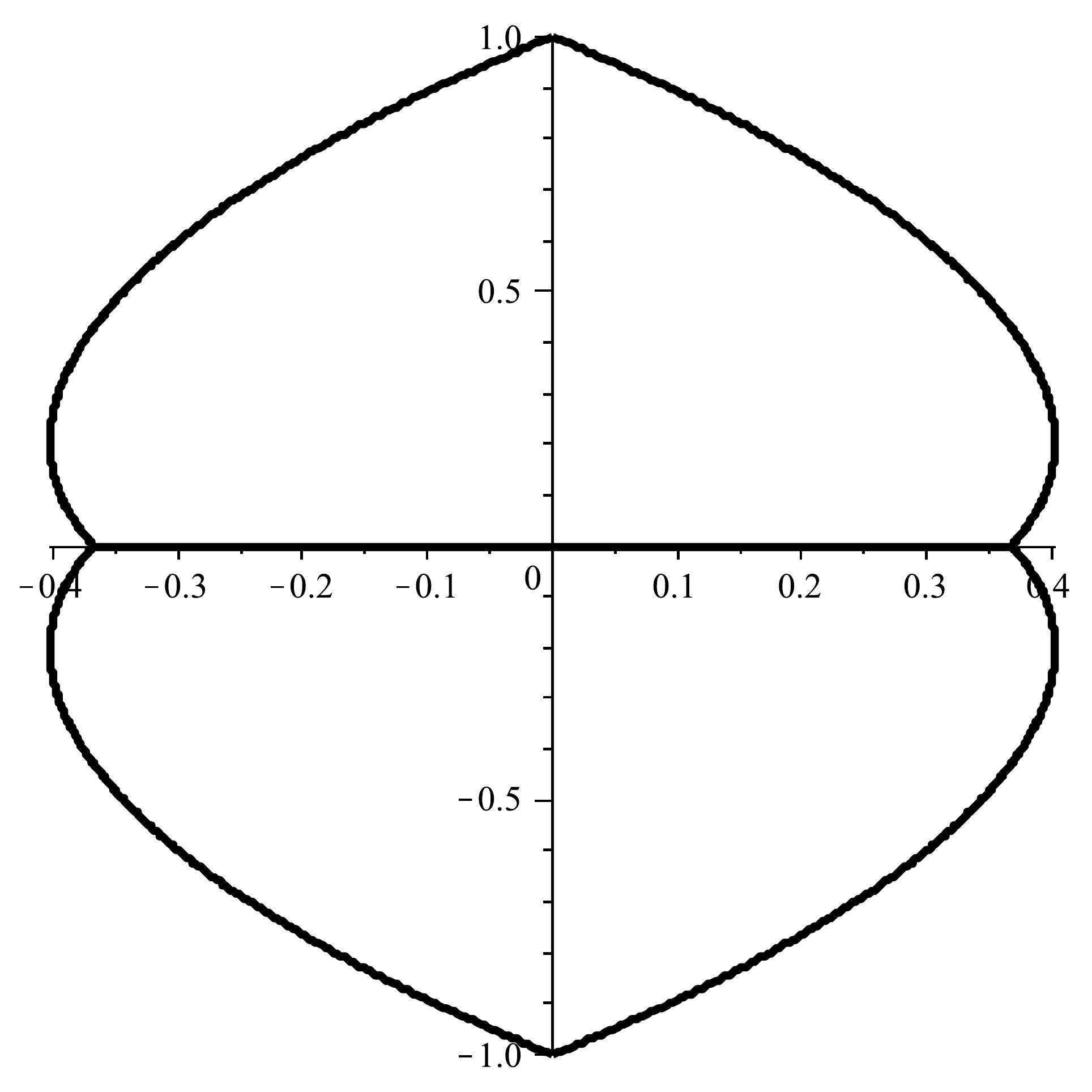} 
\caption{Zero Attractor for Taylor polynomials of $\cos(x)$ }
\label{fig:cosine}
\end{figure}

\begin{example}
{\rm
The higher order Euler polynomials $E_n^{(m)}(x)$, where $m \in {\mathbb Z}^+$, have generating
function $g(t) = (e^t+1)^m/2^m$; while the higher order Bernoulli polynomials $B_n^{(m)}(t)$
have generating function $g(t)= (e^t-1)^m / t^m$. Then their zero attractors are independent of $m$
and coincide with a scaled version of the zero attractor for the Taylor polynomials for $\cos(x)$,
see Figure \ref{fig:cosine}. 
}
\end{example}

\begin{example}
{\rm
The zero attractor for the Appell polynomials associated with generating function $g(t) =J_0(t)$, where 
$J_0(t)$ is the zero-th order Bessel function,
is a scaled version as the zero attractor for the Taylor polynomials for $\cosh(x)$, see Figure \ref{fig:bessel_attractor},
since the minimal modulus zeros of $J_0(t)$, $a=2.404825558$, are the only dominant zeros since
all the zeros of $J_0(t)$ lie on the real axis.
}
\end{example}

\begin{example}
{\rm

 Let $g(t)=(t-1) \, (t^2 +2)$. See Figure \ref{fig:appell_1} for its zero
 attractor and zeros for degree 400.

\begin{figure}[h]
\includegraphics[height=3.0in, width=3.0in]{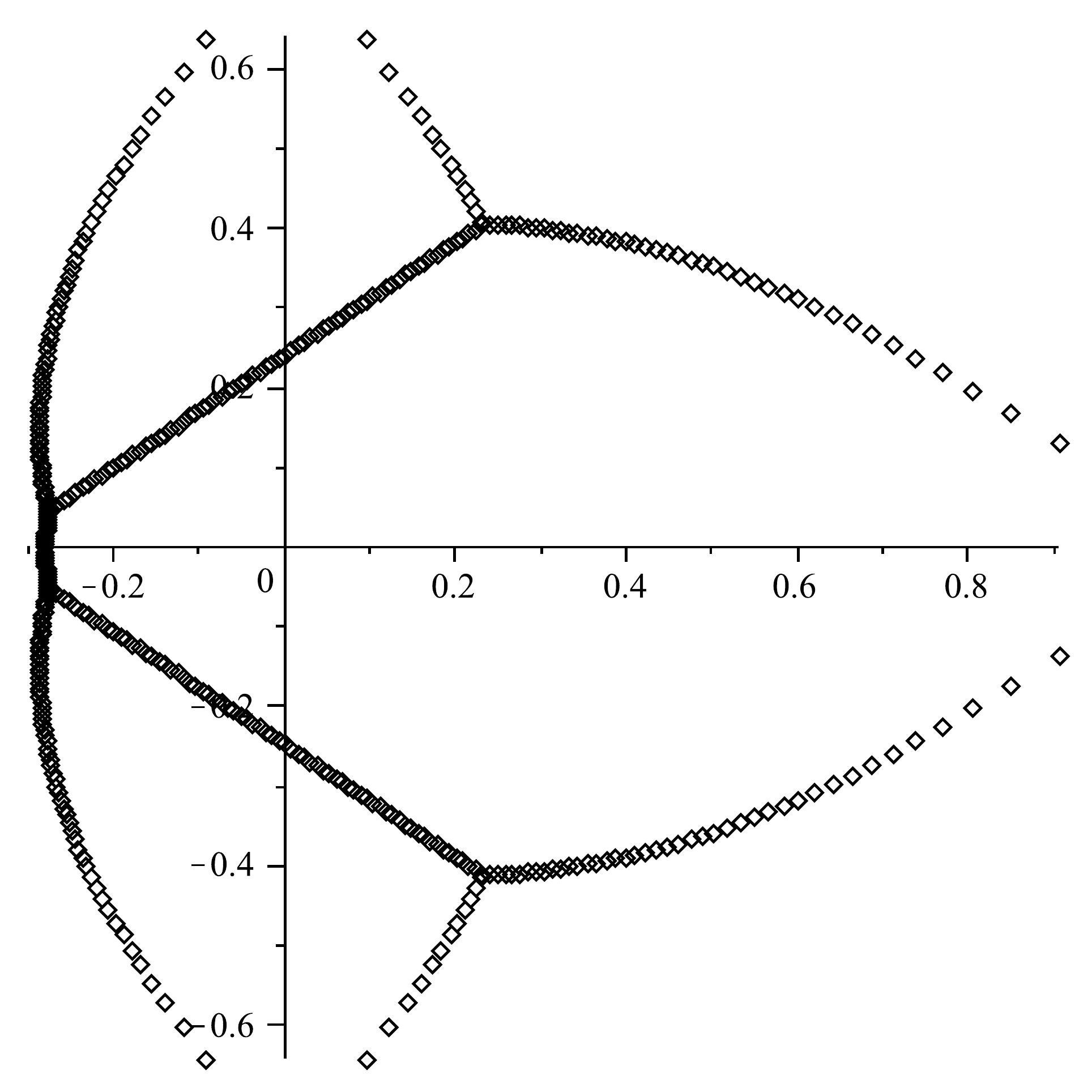} 
\quad
\includegraphics[height=3.0in, width=3.0in]{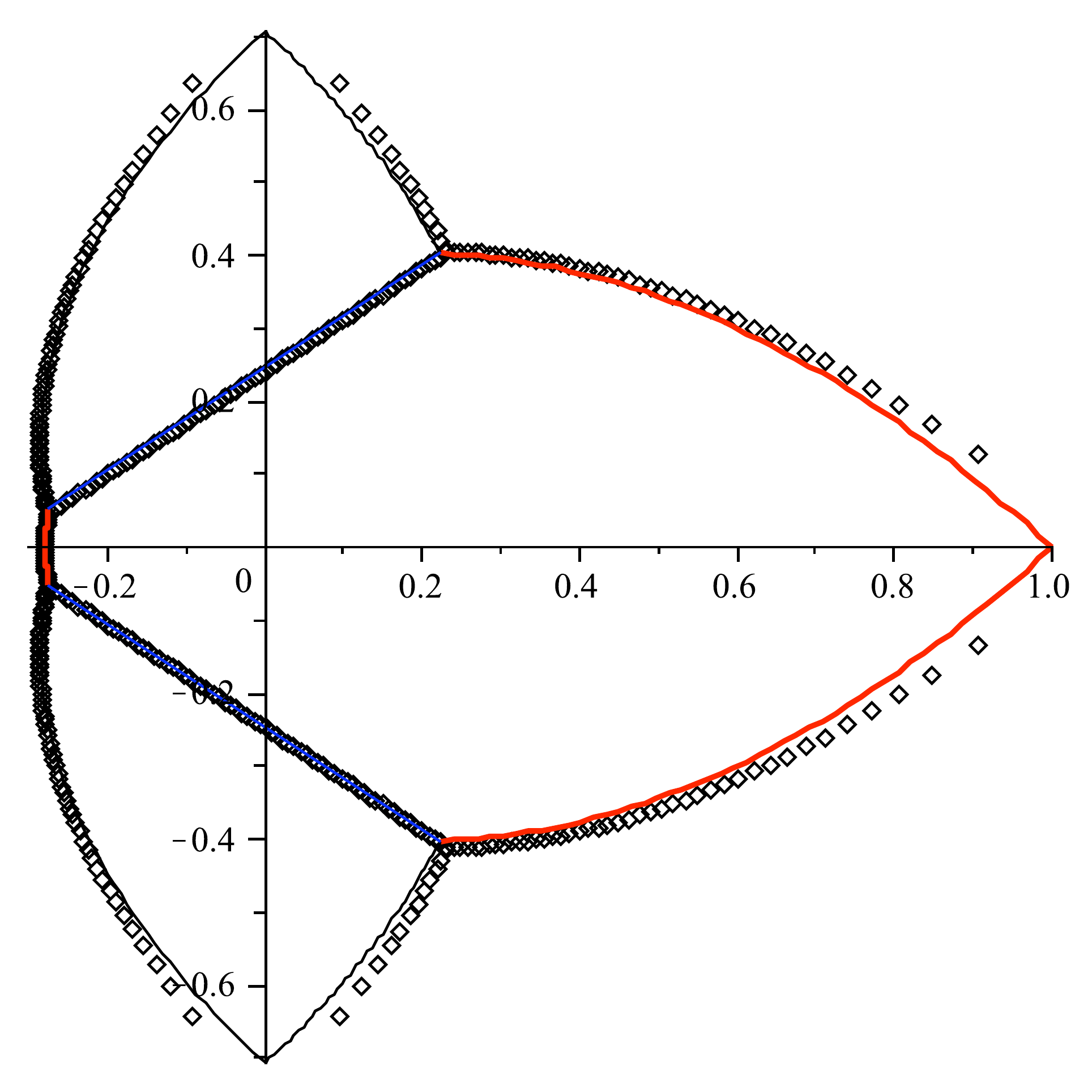} 
\caption{(a) Zeros for degree 400  polynomial with generating function $g(t)=(t-1) \, (t^2 +2)$;
(b) Zero Attractor with polynomial zeros }
\label{fig:appell_1}
\end{figure}
}
\end{example}

\begin{example}
{\rm
Consider the Appell polynomials with generating function
$g(t)= (t-a) (t-b) (t-c)$ with
$a= 1.2 e^{i 3 \pi/16}$, $b= 1.3 e^{i 7 \pi/16}$, and $c=1.5$.
 In this case, all three roots of $g(t)$ are dominant.
See Figures \ref{fig:appell_2} and \ref{fig:appell_3}.

\begin{figure} % [ht]
\includegraphics[height=3.0in, width=3.0in]{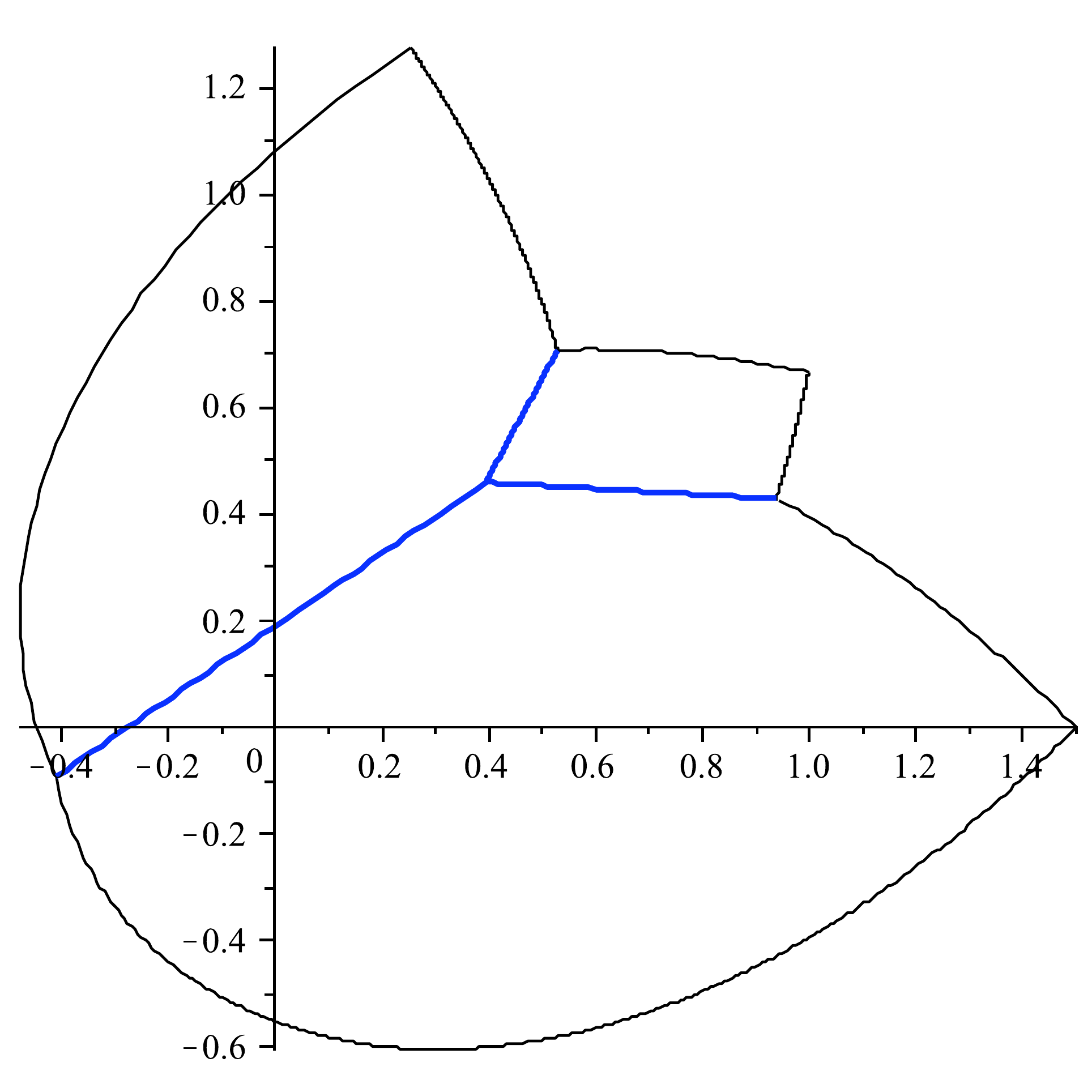} 
\quad
\includegraphics[height=3.0in, width=3.0in]{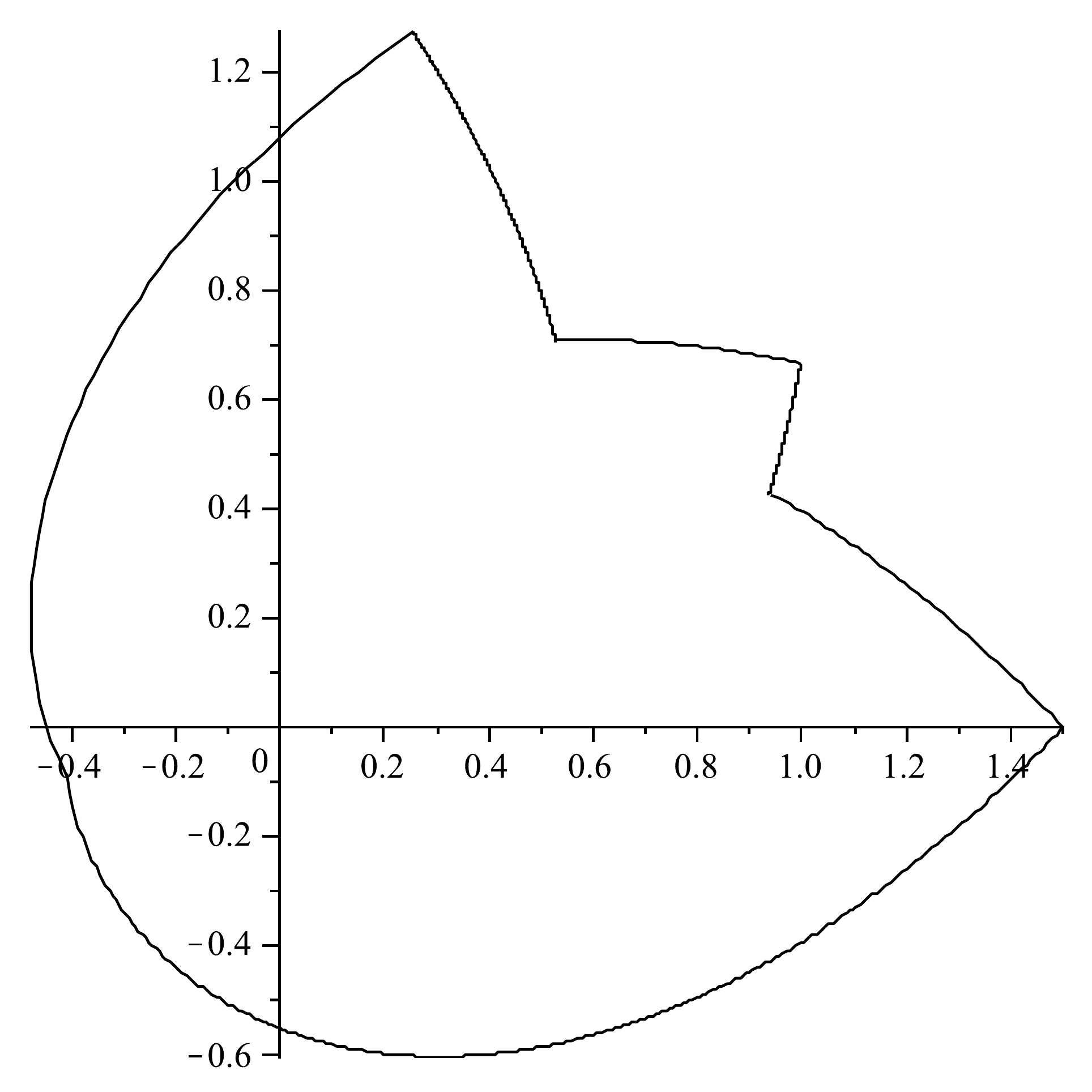} 
\caption{(a) Zero Attractor only, for generating function $g(t)=(t-a) (t-b) (t-c)$,
$a= 1.2 e^{i 3\pi/16}$, $b= 1.3 e^{i 7 \pi/16}$, $c=1.5$;
(b) Boundary of the Domain $D_0$.
 }
\label{fig:appell_2}
\end{figure}

\begin{figure} %[ht]
\includegraphics[height=3.0in, width=3.0in]{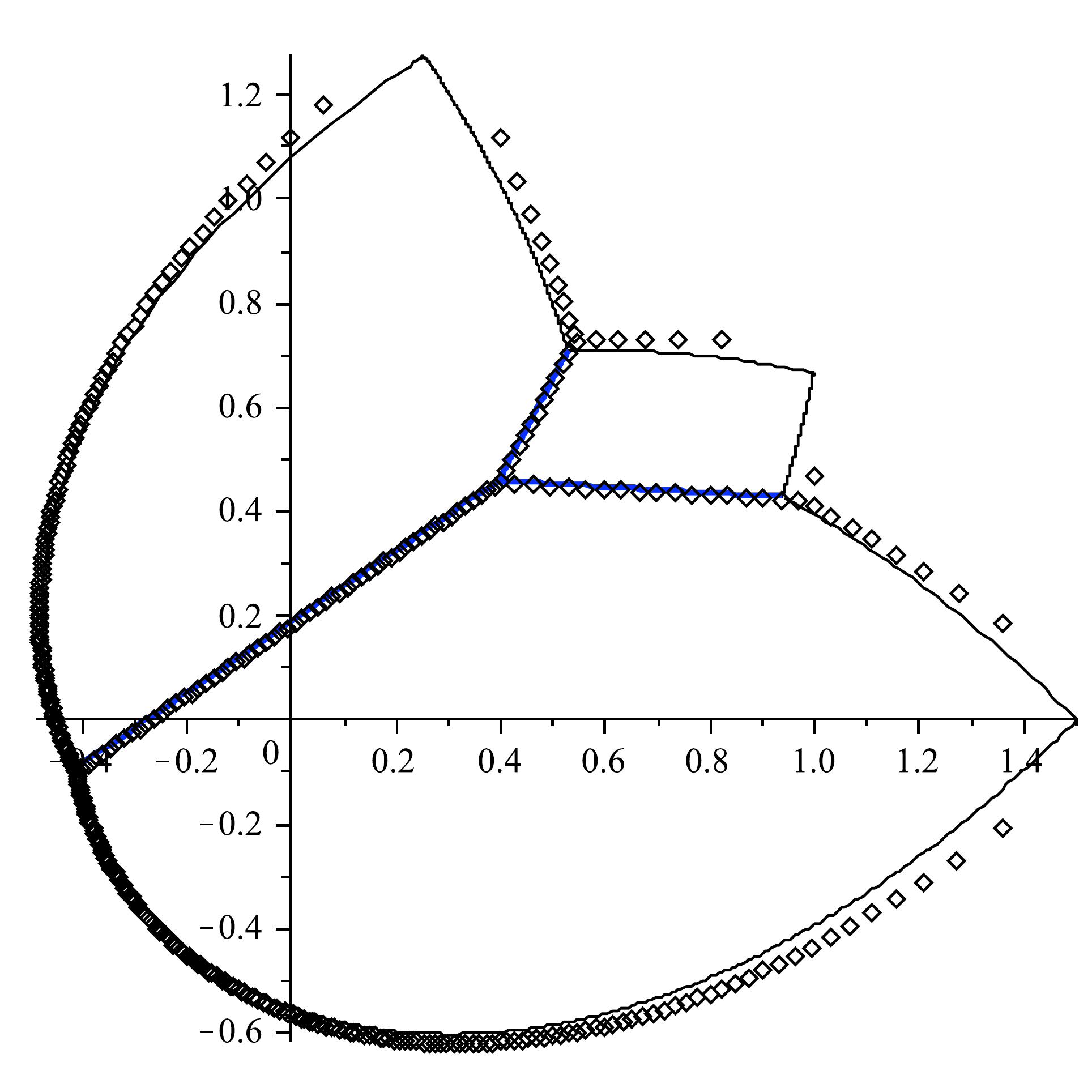} 
\caption{Zeros for degree 400 polynomial together with the Zero Attractor,
for generating function $g(t)=(t-a) (t-b) (t-c)$,
$a= 1.2 e^{i 3\pi/16}$, $b= 1.3 e^{i 7 \pi/16}$, $c=1.5$ }
\label{fig:appell_3}
\end{figure}
}
\end{example}

These last two examples both illustrate the following general fact.
We assume that the generating function $g(t)$ has exactly
three proper dominant zeros $a$, $b$, and $c$. Then
 the three
lines determined by $| \phi(ax)| = | \phi(bx)|$,
 $| \phi(ax)| = | \phi(cx)|$, and  $| \phi(bx)| = | \phi(cx)|$
 have a common intersection point, a so-called ``triple point."
 This follows by interpreting the lines as the boundary between
 the change of asymptotics of the Appell polynomial family; that is,
 the boundaries of the domains $D_a$, $D_b$, and $D_c$.

\appendix

\section{Density of Zeros}

\subsection{Introduction}

We generalize the density result for the zeros of the Euler polynomials in \cite{g-b}
to highlight how the asymptotic structure of the polynomial family may
determine the density of its zeros.

Let $\psi(x)$ be an analytic function on a domain $D \subset {\mathbb C}$
that  is conformal on $D$. 
We write $\zeta = \psi(x)$. We sometimes write $x(\zeta)$ for
$x= \psi^{-1}(\zeta)$.  

We assume that there exists $\epsilon_0>0$ and $0 \leq \alpha < \beta \leq  2 \pi$
so that the annular sector 
\begin{equation}\label{eq:def_S}  %%%%%%%%%%%%%%%%%%%%%%%%%%%%%%%  definition of region  S
S = \{ \rho e^{i \theta} : \rho \in [1- \epsilon_0,1+\epsilon_0], \theta \in [\alpha,\beta]\}
\end{equation}
lies in the image $\psi(D)$.
Next we define two subsectors of $S$ as
\begin{eqnarray*}
S_+ 
&=&
\{ \rho e^{i \theta} : \rho \in [1- \epsilon_0,1), \theta \in [\alpha,\beta]\}
\\
S_-
&=&
\{ \rho e^{i \theta} : \rho \in (1, 1+\epsilon_0], \theta \in [\alpha,\beta]\}.
\end{eqnarray*}
Let $C$ be the unimodular curve $\psi^{-1}( \{ e^{i\theta} : \theta \in [\alpha,\beta] \})$, so $| \phi(x)|=1$ for $x\in C$.
By construction,    $C$ is smoothly parametrized as
$x( e^{i \theta})$ for $\theta \in [\alpha, \beta]$.
Of course, we have
$\psi^{-1}(S) = \psi^{-1}(S_-) \cup C \cup \psi^{-1}(S_+)$ as a disjoint union.

Let $\{ T_n(x) \}$ be a sequence of analytic functions on 
$\psi^{-1}(S)$ where
we assume that the analytic functions satisfy the basic asymptotic
relation:
\begin{equation}\label{eq:T_n_asymptotics}   %%%%%%%%%%%%%%%%%%%% definition of   T_n_asymptotics
T_n(x) = 1 + a_n(x) \psi(x)^{c_n} + e_n(x),
\end{equation}
where $\{c_n\}$ is a increasing unbounded sequence of positive numbers,
 $\delta>0$ is a constant so that $|a_n(x)| \geq \delta$, and
$|a_n(x)| = \exp[ o(c_n)]$, uniformly on $\psi^{-1}(S)$.  The term $e_n(x)$ satisfy the following estimates
uniformly:
\begin{eqnarray*}
e_n(x) 
=
\left\{
\begin{array}{rl}
 o( \psi(x)^{c_n}),& \, x \in S_+,\\
 o(1),& \, x \in S_-  .
 \end{array}
 \right. 
\end{eqnarray*}
In the sequel, we may assume either form for $e_n(x)$ if $x$ lies on
the common boundary $C$ of the two regions $S_{\pm}$ that is,
$| \psi(x)|=1$.

Let $Z_n $ be the set of all zeros of $T_n$ that lie in $\psi^{-1}(S)$, which we assume is
finite for all $n$.
For $[\gamma_1,\gamma_2] \subset (\alpha, \beta)$, let
\begin{equation}\label{eq:N_n}  %%%%%%%%%%%%%%%%%%%%%%%%%%%%%   definition of   N_n
N_n( \gamma_1,\gamma_2) = \# \{ x \in Z_n : \arg x \in [\gamma_1,\gamma_2]\}.
\end{equation}
Choose $\epsilon>0$ so $3\epsilon< \epsilon_0$. By the Argument Principle,
we find that
\[
N_n(\gamma_1, \gamma_2) =
\frac{1}{2 \pi i} \int_\Gamma \frac{  \frac{d}{d \zeta} T_n(x( \zeta))}{ T_n(x(\zeta))}
\, d \zeta
\]
where $\Gamma$ is the boundary of the sector
$\{ \rho e^{i \theta} : \rho \in [1-\epsilon,1+\epsilon], \, \theta\in [\gamma_1,\gamma_2]\}$.
The closed contour 
$\Gamma$ naturally has four parts of the form
$\Gamma_{1 \pm \epsilon}$ and $\Gamma_{\gamma_j}$, $j=1,2$ where
\begin{eqnarray*}
\Gamma_{1 \pm \epsilon} &=&  \{ (1 \pm \epsilon) e^{i \theta} :  \theta \in [\gamma_1,\gamma_2]\},
\\
\Gamma_\gamma &=& \{ \rho e^{i \gamma} : \rho \in [1-\epsilon,1+\epsilon] \}.
\end{eqnarray*}

\subsection{Contributions over Arcs}

\begin{lemma}
$\displaystyle
\lim_{n \to \infty}  \int_{\Gamma_{1-\epsilon} }
\frac{  \frac{d}{d \zeta} T_n(x( \zeta))}{ T_n(x(\zeta))}
\, d \zeta =0$.
\end{lemma}

\begin{proof}
Let $\zeta \in \Gamma_{1- \epsilon}$. Then $x(\zeta) \in \psi^{-1}(S_-)$
and $e_n(x(\zeta)) \to 0$ uniformly on  $\Gamma_{1- \epsilon}$.
Furthermore, we find
\[
| a_n ( x(\zeta)) \psi^n( x( \zeta))|
\leq \exp[ o(c_n)]  | \zeta|^{c_n} =
e^{ - c_n \ln( 1/(1-\epsilon))+ o(c_n))} = o(1).
\]
From the basic asymptotic expression
\[
T_n( x( \zeta)) = 1+ a_n( x(\zeta)) \zeta^{c_n} + e_n ( x(\zeta)),
\]
we find that $\ln[ T_n( x (\zeta))]$ converges uniformly to 0 on $\Gamma_{1-\epsilon}$.
Hence their derivatives must also converge uniformly to 0 and so the desired integrals
converge to 0.
\end{proof}

\begin{lemma}
$\displaystyle
\lim_{n \to \infty} \Im
\left(
  \frac{1}{c_n} \int_{\Gamma_{1+\epsilon} }
\frac{  \frac{d}{d \zeta} T_n(x( \zeta))}{ T_n(x(\zeta))}
\, d \zeta
\right)
 = \gamma_2-\gamma_1$.
\end{lemma}

\begin{proof}
For $\zeta \in \Gamma_{1+\epsilon}$, $| \zeta| = 1+\epsilon$ and
$x(\zeta) \in \psi^{-1}(S_+)$. By the basic expansion
\[
T_n( x( \zeta)) = 1+ a_n( x(\zeta)) \zeta^{c_n} + e_n ( x(\zeta)),
\]
we have that
\[
\frac{ T_n( x(\zeta))}{ a_n( x(\zeta)) \zeta^{c_n}}
=
1 + \frac{ 1+e_n( x(\zeta))}{ a_n(x(\zeta)) \zeta^{c_n}}
=
1+ \frac{  \zeta^{-c_n} + \frac{ e_n(x(\zeta)) } {\zeta^{c_n}}} {a_n(x(\zeta))}.
\]
We recall our assumptions 
that $ | a_n( x( \zeta))| \geq \delta$ and $ e_n( x(\zeta))/ \zeta^{c_n} = o(1)$. 
Hence, we find that
\[
\left|  \frac{ 1+e_n( x(\zeta))}{ a_n(x(\zeta)) \zeta^{c_n}} \right|
=
\left|    \frac{  \zeta^{-c_n} + \frac{ e_n(x(\zeta)) } {\zeta^{c_n}}} {a_n(x(\zeta))}    \right|
\leq
\left|
\frac{ (1+\epsilon)^{c_n} + o(1)}{ \delta}
\right| = o(1)
\]
uniformly on $\Gamma_{1+\epsilon}$.
In particular, $\frac{T_n( x(\zeta)) }{ a_n(x(\zeta)) \zeta^{c_n}}$ converges to $1$
uniformly on $\Gamma_{1+\epsilon}$ so 
$\ln[ \frac{T_n( x(\zeta)) }{ a_n(x(\zeta)) \zeta^{c_n}}]$ 
converges uniformly to 0 there as do their derivatives. In other words,
we know that
\[
\frac{  \frac{d}{d \zeta} T_n(x( \zeta))}{ T_n(x(\zeta))} 
-
\frac{  \frac{d}{d \zeta} a_n( x(\zeta))}{ a_n(x(\zeta))} - \frac{c_n}{\zeta} \to 0
\]
Note that $| a_n(x)| = \exp[ o(c_n)]$ implies that
$| a_n'(x)/a_n(x)| = o(c_n)$.
Hence, we find
\begin{eqnarray*}
\frac{  \frac{d}{d \zeta} T_n(x( \zeta))}{ T_n(x(\zeta))} 
&=&
\frac{  \frac{d}{d \zeta} a_n( x(\zeta))}{ a_n(x(\zeta))} + \frac{c_n}{\zeta}+o(1)
\\
&=&
\frac{  \frac{d}{dx} a_n( x)}{ a_n(x)}\, \frac{dx}{d \zeta}
  + \frac{c_n}{\zeta} + o(1)
  \\
  &=&
  \frac{c_n}{\zeta} + o(c_n).
\end{eqnarray*}
We now conclude that
\[
\frac{1}{c_n} 
\frac{  \frac{d}{d \zeta} T_n(x( \zeta))}{ T_n(x(\zeta))} \to \frac{1}{\zeta}
\textrm{  uniformly on   }  \Gamma_{1+\epsilon}.
\]
The lemma now follows easily.
\end{proof}

\subsection{Backlund's Method}

Our estimates for the integrals over $\Gamma_\gamma$,
for $\gamma \in [\gamma_1,\gamma_2]$, are inspired by
the 1918 method of R. Backlund's proof of the Riemann-von Mangoldt
asymptotic  formula 
for the number of zeros of the Riemann zeta function. We follow the exposition
of Chandrasekharan \cite{chandra} (pages 35-38).

Recall that $\Gamma_\gamma$ is parametrized as
$\rho e^{ i \gamma}$ for $\rho \in [1-\epsilon,1+\epsilon]$.

Let $\ell$ be the number of zeros of $\Re[ T_n( x(\zeta))]$ for $\zeta \in \Gamma_\gamma$ exclusive of endpoints. 
Then the contour integral
along $\Gamma_\gamma$ can be written as a sum of integrals over line
segments
$C_{ab} $ paramerized as $\rho e^{i \gamma}$, $\rho \in [a,b]$ where
 $a$ and $b$ are two consecutive zeros of $\Re[ T_n( x(\zeta))]$.
In particular,  $\Re[ T_n( x( \rho e^{i \theta} ))]$ has constant sign for $\rho \in [a,b]$.
Then 
\[
\Im
\left(
\int_{C_{ab}} 
\frac{  \frac{d}{d \zeta} T_n(x( \zeta))}{ T_n(x(\zeta))}
\, d \zeta 
\right)
=
\Im \left( \int_{C_n} \frac{d \xi}{\xi} \right),
\]
where $C_n$ is the image of the line segment $\rho e^{i \gamma}$,
for $\rho \in [a,b]$, under the map
$\zeta \mapsto T_n(x(\zeta))$. 
By assumption, $C_n$ can only intersect
the imaginary axis $i {\mathbb R}$ only at its two endpoints; in particular,
$C_n$ must lie either in the left and right half-plane. By Cauchy's theorem,
we can deform $C_n$ into a semicircle $K_{ab}$ that lies in the same half-plane
and has the same endpoints on $i {\mathbb R}$ so that the value of the contour
integral is unchanged. This allows us to make the estimate
\[
\left| \Im \int_{C_n} \frac{d \xi}{ \xi} \right| =
\left|\Im  \int_{K_{ab}} \frac{d \xi}{ \xi}  \right| \leq \pi.
\]
We summarize this discussion as:
\begin{lemma}
$\displaystyle
\left| \Im \int_{\Gamma_\gamma}  
\frac{  \frac{d}{d \zeta} T_n(x( \zeta))}{ T_n(x(\zeta))}
\, d \zeta 
 \right| \leq ( \ell+1) \pi
$ where $\ell$ is the number of zeros of 
$\Re[ T_n( x( \rho e^{i \theta} ))]$ for $\rho \in [1-\epsilon,1+\epsilon]$.
\end{lemma}

We now use Jensen's formula to make a useful estimate for $\ell$.

For an analytic function $h(z)$ on some domain $E$, define 
$\tilde{h}(z)$ on $E_c$ as $\tilde{h}(z)= \overline{ h( \overline{z})}$
where $E_c = \{ \overline{z} : z \in E\}$ which will be analytic on $E_c$.

For $\xi \in D(1-\epsilon; 2 \epsilon)$, which is symmetric about the real axis,
let
\[
\hat{T}_n(\xi)
=
\tfrac{1}{2}
\left[
T_n( x( \xi e^{i \gamma})) + \tilde{T}_n( \tilde{x}( \xi e^{- i \gamma}))
\right].
\]
Then $\hat{T}_n(\xi)$ is analytic on $D(1-\epsilon; 2 \epsilon)$
and 
\[
\hat{T}_n(\xi) = \Re[ T_n( x( \xi e^{i \gamma}))], \quad
\xi \in [1-\epsilon,1+\epsilon] \subset {\mathbb R}.
\]
For convenience, we recall Jensen's inequality.
Let $h(z)$ be an analytic function on the closed disk $\overline{D}(a;R)$,
and let $0<r<R$. 
Suppose $h(a) \neq 0$.
Let $m$ be the number of zeros of $h(z)$ in the
closed disk $\overline{D}(a;r)$ counted according to their multiplicity.
Then
\[
\left( \frac{R}{r}\right)^m
\leq
\frac{  \max \{ | h(z) | : |z-a|=R \} } { | h(a)|}.
\]
Now each zero of $\Re[ T_n( x(\zeta))]$ for $\zeta \in \Gamma_\gamma$
corresponds to a zero of $\hat{T}_n( \xi)$ for $\xi \in [1-\epsilon,1+\epsilon]$.
Let $\hat{\ell}$ be the number of zeros of $\hat{T}_n(\xi)$ 
for $\xi \in D(1-\epsilon; 2 \epsilon)$. Then we have at once the inequality
\[
\ell \leq \hat{\ell}.
\]
We will apply Jensen's inequality to the disk $D(1-\epsilon; 3 \epsilon)$
and $r = 2 \epsilon$ to obtain
\[
\left( \frac{3}{2} \right)^{ \hat{\ell}} 
\leq 
\frac{  \max \{ |  \hat{T}_n(\xi) | : |\xi-(1-\epsilon)|=3 \epsilon \} } { |  \hat{T}_n(1-\epsilon )|}.
\]

\begin{lemma}
(a) $\hat{T}_n(1-\epsilon) = 1+o(1)$.
\\
(b) 
$\displaystyle \max \{  | \hat{T}_n( \xi) :  | \xi - (1-\epsilon)|=3 \epsilon  \}  
= O\left( e^{ o(c_n)} (1+ 2 \epsilon)^{c_n} \right).$
\end{lemma}

\begin{proof}
(a)
Since $(1-\epsilon) e^{i \gamma} \in S_-$, we have the estimate
\[
|  T_n( x( (1-\epsilon) e^{i \gamma}))| \leq 1+ e^{ o(c_n)}(1-\epsilon)^{c_n} + o(1)
= 1 +o(1).
\]
In particular, $|\hat{T}_n( 1-\epsilon) | = 1+ o(1)$.

(b)
To estimate the maximum of $|\hat{T}_n( \xi)|$ for $| \xi - (1-\epsilon)|=3 \epsilon$,
we observe that
\[
| \hat{T}_n( \xi)| \leq
\tfrac{1}{2} 
\left[
\left| T_n(  x(  \xi e^{i \gamma})) \right| 
+ 
\left| \tilde{T}_n( \tilde{x} ( \xi e^{- i \gamma})) \right|
\right].
\]
Let  $\zeta$ lie in the closed disk $\overline{D}(1-\epsilon;3\epsilon)$
so $| \zeta| \leq 1+ 2\epsilon$.
By the basic asymptotic relation, we find that whether $\zeta \in S_+$ or
$S_-$:
\begin{eqnarray*}
| T_n(  x( \zeta  ))  | 
&\leq&
| 1+a_n(x(\zeta)) \zeta^{c_n} + e_n(x(\zeta))|
\\
&\leq&
1+ e^{ o(c_n)} (1+ 2 \epsilon)^{c_n} + o( (1+2 \epsilon)^{c_n})
= O \left( e^{o(c_n)} \,(1+ 2 \epsilon)^{c_n} \right),
\end{eqnarray*}
where the big-oh constant holds uniformly on  
$\overline{D}(1-\epsilon;3\epsilon)$.
In particular, this estimate holds for  $\xi$ that lie on the circle
$| \xi - (1-\epsilon)|=3 \epsilon$.
A similar estimate holds for $| T_n(  \tilde{x}( \xi e^{i \gamma}))|$.
We sum up this discussion as
\[
\max \{  | \hat{T}_n( \xi) :  | \xi - (1-\epsilon)|=3 \epsilon  \}  
= O\left( e^{ o(c_n)} (1+ 2 \epsilon)^{c_n} \right).
\]
\end{proof}

The last lemma together with Jensen's inequality
 allows us to make an estimate for $\hat{\ell}$:
\[
\left( \frac{ 3}{2} \right)^{ \hat{\ell}}
\leq 
\frac{   \max \{  | \hat{T}_n( \xi) :  | \xi - (1-\epsilon)|=3 \epsilon  \}    }
{  \hat{T}_n( 1-\epsilon)}
=
O\left( e^{ o(c_n)} (1+ 2 \epsilon)^{c_n} \right).
\]
Recalling that $\ell \leq \hat{\ell}$, we have the bound
\[
\ell \leq \frac{ \ln[   O\left( e^{ o(c_n)} (1+ 2 \epsilon)^{c_n} \right)   ] }{ \ln(3/2)}.
\]
Hence, we have the estimate for the integral
\[
\frac{1}{c_n}
\left|
\Im \int_{\Gamma_\gamma}  
\frac{  \frac{d}{d \zeta} T_n(x( \zeta))}{ T_n(x(\zeta))}
\, d \zeta 
\right|
\leq \frac{\pi}{\ln(3/2)}  \ln( 1+ 2 \epsilon) + O\left( \frac{1}{c_n}\right).
\]
Since these bounds hold for all $\epsilon>0$ sufficiently small, we have
shown the
\begin{lemma}
For all $\epsilon>0$ sufficiently small,
\[
\frac{1}{c_n}
\left| \Im
\int_{\Gamma_\gamma} 
\frac{  \frac{d}{d \zeta} T_n(x( \zeta))}{ T_n(x(\zeta))}
\, d \zeta
\right| \leq  \frac{\pi}{\ln(3/2)}  \ln(1+2 \epsilon) + O \left( \frac{1}{c_n}  \right)
\]
where $\Gamma_\gamma$ is the line segment
$\rho e^{i \gamma}$, $\rho \in [1-\epsilon,1+\epsilon]$.
\end{lemma}

Finally, this last inequality allows us to make the estimates
\[
\frac{\gamma_2-\gamma_1}{2 \pi} -   \frac{1}{ 2\ln(3/2)} \ln( 1+ 2 \epsilon) + 
O\left( \frac{1}{c_n}\right)
\leq \frac{ N_n( \gamma_1,\gamma_2)}{ 2\pi c_n} 
\leq
\frac{\gamma_2-\gamma_1}{2 \pi} +  \frac{1}{ 2\ln(3/2)} \ln( 1+ 2 \epsilon)+ 
O\left( \frac{1}{c_n}\right).
\]

\subsection{Main Density Theorem}

By combining the above lemmas and noting that
these bounds hold for all $\epsilon>0$ sufficiently small, we obtain our main density result.

\begin{theorem}\label{thm:density}
Let $\alpha < \gamma_1 < \gamma_2 < \beta$, and let
$N_n(\gamma_1,\gamma_2)$ denote the number of zeros of $T_n(x)$ whose arguments lie in $[\gamma_1,\gamma_2]$,
given in equation  (\ref{eq:N_n}). Then
\[
\lim_{n \to \infty} \frac{ N_n(\gamma_1,\gamma_2)}{c_n} 
=
\frac{ \gamma_2-\gamma_1}{2 \pi};
\]
that is, the image of the zero density under $\psi$ is Lebesgue measure on an arc of the
unit circle.
\end{theorem}

We need to recall the notions of $\limsup$ and $\liminf$ of a sequence $\{X_n\}$ of compact
sets in the complex plane. Now $x^* \in \limsup X_n$ if for every neighborhood $U$ of $x$,
there exists a sequence $x_{n_k} \in X_{n_k} \cap U$  that converges to $x^*$ while $x^* \in
\liminf X_n$ if for every neighborhood $U$ of $x$, there exists an index $n^*$ and a sequence
$x_n \in X_n \cap U$, for $n \geq n^*$ that converges to $x^*$.

It is known that
if the $\liminf X_n$ and $\limsup X_n$ agree and are uniformly
bounded, then the sequence $\{X_n \}$ converges in the Hausdorff metric.

When the density result holds, then the $\liminf Z(T_n)$ must agree with
$\limsup Z(T_n)$. 
Hence, we have the following:

\begin{corollary}
As compact subsets of $\psi^{-1}(S)$, $Z(T_n)$ converges to the unimodular curve $C$ in the Hausdorff metric.
\end{corollary}

Although we can determine the zero attractor and the zero density completely
in the above framework. it is conceptually useful to have the result of Sokal
that gives a description of the support of the zero density measure.

{\sl
{\bf [Sokal] \cite{sokal}:} \label{thm:sokal}
Let $D$ be a domain in ${\mathbb C}$, and let $z_0 \in D$.
Let $\{ g_n\}$ be analytic functions on $D$,  and let
$\{a_n\}$ be positive real constants such that
$\{ |g_n|^{a_n} \}$ are uniformly bounded on the compact
subsets of $D$. Suppose that there does {\rm not}
exist a neighborhood $V$  of $z_0$ and a function
$v$ on $V$ that is either harmonic or else identically
$- \infty$ such that
$\displaystyle \liminf_{n \to \infty} a_n \ln |g_n(z)| \leq v(z) \leq
\limsup_{n \to \infty} a_n \ln |g_n(z)|$ for all $z \in V$.
Then $z_0 \in \liminf {Z}(g_n)$.
}

{\sl Remark}:
We can state the asymptotic form for $T_n(x)$ in a more symmetric form as:
\[
T_n(x) = \psi_0(x) + \sum_{k=1}^N a_{n,k}(x) \psi_k(x)^{c_n}
+ e_n(x)
\]
where $N$ is fixed and the error term has the form
\[
e_n(x) = o( \max \{\psi_k(x)^{c_n},  0 \leq k \leq n \} )
\]
This version explains the asymmetry in our first result
where we have $\psi_0(x) =1$ and  the zeros
accumulate along the curve $| \psi_0(x)| = |\psi(x)|$.

\subsection{Special Case}

Theorem \ref{thm:density} shows that the images of zeros 
under the conformal map $\psi$
are uniformly distributed along the corresponding circular arc. 
 It can be
applied to many cases that arise in a broad spectrum. A special case is
worthy of attention; namely, the analytic arc $C$ is a straight line segment
and $\psi(x)$
 (see Lemma \ref{lemma:5.1})
has  the form $e^{ax+b}$, where $a$ and $b$ are
constants. 

\begin{corollary}
\label{cor:constant} If the analytic arc $C$ is a straight line segment and 
$\psi(x) $ is of the form $e^{ax+b}$, where $a$ and $b$ are constants, then the
zero density along the line segment $C$ is a multiple of Lebesgue measure. 
\end{corollary}

\begin{proof}  % proof 25
 Let $c_1$ and $c_2$ be the endpoints of the line segment so it is parametrized
 as $x(t)=c_1(1-t) + c_2 t$, $0\leq t \leq 1$. 
 Then $\psi(x(t))$ becomes
$
\psi(x(t))= e^{ a [ c_1 (1-t) + c_2 t] + b}
$
Since $| \psi( x(t))|=1$, $\psi( x(t))$ can be written as
\[
\psi( x(t)) = e^{ 2 \pi i \theta(t)},
\]
where $\theta(t)$ is a linear function of $t$.
By Theorem \ref{thm:density} the density
of images of zeros under $\psi(x)$ along the corresponding circular arc is 
Lebesgue measure. Hence, its pull-back under $\psi$ is also Lebesgue measure
since $\theta(t)$ is linear.
\end{proof}

\end{document}